\documentclass[12pt,a4paper,oneside]{amsart}
\usepackage{fullpage,url,amssymb,enumerate,colonequals,array}
\newcolumntype{P}[1]{>{\raggedright\arraybackslash}p{#1}} 
\usepackage[all,pdf]{xy} 

\usepackage{comment}
\usepackage{color}
\usepackage{charter,euler} 
\DeclareSymbolFont{bchoperators}{T1}{bch}{m}{n}
\SetSymbolFont{bchoperators}{bold}{T1}{bch}{b}{n}
\makeatletter
\renewcommand{\operator@font}{\mathgroup\symbchoperators}
\makeatother

\DeclareFontEncoding{OT2}{}{} 

\usepackage{titlesec} 
\titleformat{\section}{\normalfont\bfseries\filcenter}{\thesection}{1em}{}
\titleformat{\subsection}{\normalfont\bfseries}{\thesubsection}{1em}{}
\titleformat{\subsubsection}{\normalfont\itshape}{\thesubsubsection}{1em}{}

\definecolor{darkgreen}{rgb}{0,0.5,0}

\usepackage[
        colorlinks, citecolor=darkgreen,
        pdfauthor={Michael Stoll},
        pdftitle={Documentation for the ratpoints program},
]{hyperref}

\newcommand{\Q}{{\mathbb Q}}
\newcommand{\R}{{\mathbb R}}
\newcommand{\BP}{{\mathbb P}}

\newcommand{\rpversion}{3.1.0} 

\begin{document}

\title{Documentation for the ratpoints program}

\author{Michael Stoll}
\address{Mathematisches Institut,
         Universit\"at Bayreuth,
         95440 Bayreuth, Germany.}
\email{Michael.Stoll@uni-bayreuth.de}
\date{September 24, 2026}

\maketitle


\section{Introduction}

This paper describes the \texttt{ratpoints} program. This program tries to
find all rational points within a given height bound on a hyperelliptic
curve in the most efficient way possible.


\section{History and Acknowledgments}

This program goes back to an implementation of the `quadratic sieving'
idea by \textbf{Noam Elkies} that was around in the early 1990s. My own
first contribution was to replace the \texttt{char} arrays that were used
to store the sieving information by bit arrays, in~1995. \textbf{Colin Stahlke}
then made use of the \texttt{gmp} library, so that points could be checked
exactly, and implemented the selection of sieving primes according to
their likely success rate, in~1998. After that, I successively put in
numerous improvements (and some bug fixes). For details of how the
program works, see Section~\ref{impl} below.

Along with \textbf{Noam Elkies} and \textbf{Colin Stahlke}, I would like to
thank \textbf{John Cremona} and \textbf{Sophie Labour} for bug reports and
suggestions for improvements. Further thanks are due to \textbf{Bill Allombert}
for a first implementation of the use of 256-bit (and possibly 512-bit) registers.

From September~2026, I have been using \textbf{Claude Code} (in various
incarnations: Opus~5, Fable~5, Fable~5.1) to optimize the code and find and fix
bugs. This has resulted in version~3.0.0, which includes
a whole array of optimizations compared to the 2.2.x~versions, making
the program run 2 to~10 times faster (depending on the kind
of curve and the height bound). Because of this substantial overall
improvement and the change in the development workflow, I decided to
increase the major version number from 2 to~3, even though many of
the optimizations are incremental and there are few really new ideas involved.
The current version~3.1.0 adds the option to run the sieve on several
threads (Sections~\ref{sec:threads} and~\ref{sec:library}).

Note that version~2.2.4 is the last one that supports 32-bit architectures;
versions 3.0.0 and up need a word size of 64~bit.


\section{Availability}

The \texttt{ratpoints-\rpversion} package can be downloaded from my homepage,
see~\cite{ratpoints}, or obtained from \texttt{github} at
\url{https://github.com/MichaelStollBayreuth/ratpoints}.

This program is free software: you can redistribute it and/or modify
it under the terms of the GNU General Public License as published
by the Free Software Foundation, either version~2 of the License, or
(at your option) any later version.

This program is distributed in the hope that it will be useful,
but \textbf{without any warranty}; without even the implied warranty of
\textbf{merchantability} or \textbf{fitness for a particular purpose}. See the
GNU General Public License for more details.

You should receive a copy of the GNU General Public License
along with this program.

If not, see \url{http://www.gnu.org/licenses/}.


\section{Installation}

This section describes the installation procedure under Linux.

\subsection{Extract the archive}

\begin{verbatim}
> tar xzf ratpoints-3.1.0.tar.gz
\end{verbatim}

This sets up a directory \texttt{ratpoints-\rpversion} containing the various
files that belong to the installation. Alternatively, clone the github
repository mentioned above into a directory of your choice (and replace
\texttt{ratpoints-\rpversion} by the name of that directory in the instructions below).

\subsection{Build the program and library}\label{sec:build}

Do

\begin{verbatim}
> cd ratpoints-3.1.0
> make all
\end{verbatim}

This will build the library \texttt{libratpoints.a}, the executable
\texttt{ratpoints}, and also run pdf\/latex to generate this documentation.

Five of the constants that govern the sieve depend on the machine rather than
on the curve. Reasonable values are compiled in, but \texttt{make tune} will
measure better ones for your machine if you want it to; see
Section~\ref{sec:retune}. It is a separate step because it takes several
minutes.

The sieve can run on several threads (option~\verb+-t+,
Section~\ref{sec:threads}; field \texttt{num\_threads},
Section~\ref{sec:library}), for which the library uses POSIX threads and
is compiled with \verb+-pthread+. On a system without them,
\begin{verbatim}
> make THREADS=0
\end{verbatim}
builds the library without threads (\verb+-DRATPOINTS_NO_THREADS+); the
option and the field are then accepted and ignored, and the sieve runs on
the calling thread. At most $33$ threads sieve at once,
whatever is asked for: the correction of the sieving parameters during a
run (Section~\ref{sec:adapt}) takes effect $32$ blocks of the loop over
the denominators after it is decided, which is what lets that many blocks
be in progress at the same time, and the library creates no threads
beyond that. On a machine with more cores, build with a larger lag, as in
\begin{verbatim}
> make all CCFLAGS=-DRP_ADAPT_LAG=64
\end{verbatim}
for up to $65$ threads (the \texttt{Makefile} notices the changed flags
and rebuilds what depends on them). The correction then comes later in the run, by a
few per cent of it; the points found are the same, but the reports of
\verb+-v+ on a long run may differ from those of a build with the
default lag, as they may between any two builds.

By default, the program has the $53$ odd primes below~$256$ available and
considers the smallest~$30$ of them; for a curve that needs more, it will
look further of its own accord, as described in Section~\ref{sec:spchoice}.
If you intend to do computations with curves that have many rational points,
it may make sense to make more primes available. This can be achieved via

\begin{verbatim}
> make all PRIME_SIZE=9
\end{verbatim}

The \texttt{PRIME\_SIZE} argument can be given any value from 5 to~10;
the default value is~8. The precise meaning is that the program will
work with primes $< 2^s$, when $\text{\texttt{PRIME\_SIZE}} = s$. The
\texttt{Makefile} notices that the setting has changed and rebuilds what
depends on it. Only the primes a run actually asks for get sieving tables,
so a larger table costs nothing until it is used (Section~\ref{fd}).

The sieve works on several numerators at a time, in a register whose width is
selected by the line
\begin{verbatim}
CCFLAGS1 = ${CCFLAGS256}
\end{verbatim}
in the \texttt{Makefile}. It can be set to \verb+${CCFLAGS64}+ (64-bit
registers, i.e., plain \texttt{unsigned long}s, which works on any machine
the library builds on at all), \verb+${CCFLAGS128}+,
\verb+${CCFLAGS256}+ (the default, using 256-bit AVX registers) or
\verb+${CCFLAGS512}+. The \texttt{Makefile} notices that the setting has
changed and rebuilds everything that depends on it, so a \texttt{make
distclean} in between is not necessary. Wider
registers help the first stage of the sieve, where most of the time of a
long run goes, but they help the later stages much less and the set-up not
at all, so the gain is nowhere near proportional to the width: on the
author's machine, the random test curves at a height bound of $2\cdot10^5$
take $1.24$ units of time with 128~bits for every unit with 256~bits. The
plain 64-bit build takes $3.7$ units; it sieves one word at a time and leaves
the rest to the compiler, and is meant for machines that offer nothing
wider.


You may have to change \verb+-mavx2+ in the line
\begin{verbatim}
CCFLAGS256 = -DUSE_AVX -mavx2
\end{verbatim}
into \verb+-mavx+ or leave it out altogether, depending on the capabilities
of your CPU. You will notice that AVX2 (or AVX) is not supported, when the
executable throws an `unknown instruction' error or similar.

There is also a variant using 512-bit registers, activated by
\begin{verbatim}
CCFLAGS1 = ${CCFLAGS512}
\end{verbatim}
which requires a CPU with AVX512F capability. I have no such CPU available
and can test this variant only in the emulation described below;
\textbf{Drew Sutherland} has run the complete test suite of version~3.0.0 on
an AMD Ryzen~9 9950X3D (Zen~5 architecture), where it passes. It is worth
timing the 256-bit and the 512-bit version against each other on your
machine: the sieving loop is limited by how fast it loads bit arrays from the
first-level cache and by how much of that cache the sieve tables occupy,
rather than by arithmetic, so wider registers do not automatically pay off,
and some Intel server processors reduce their clock frequency when 512-bit
instructions are used. On the Zen~5 machine, with the constants of
\texttt{make tune} (Section~\ref{sec:retune}) in both builds, the 512-bit
version was 14\% faster than the 256-bit one on a long run (the curve of
\texttt{make timing}, height bound $4\cdot10^5$), equal on the random test
curves at~$16383$ and on the sextic of \texttt{make test2} at~$10^6$, and
8\% slower on the point-rich curves at~$16383$, where building the sieve
tables, twice as large with 512-bit registers, is a large share of the run.
Wider registers thus help long runs on curves that need few sieving moduli,
and not short ones. Run \texttt{make tune} after choosing the width: on that
machine it gained 6 to 7\% at either width over the constants compiled in,
which were measured on the author's laptop, and it found the same values at
both widths.


\subsection{Run various tests}\label{sec:tests}

Run

\begin{verbatim}
> make test1
\end{verbatim}

in the working directory. This will build an executable \texttt{rptest} and
then run (and time) it. Finally, the output (which was written to a file
\texttt{rptest.out}) is compared against \texttt{testbase}, which contains the
output of a sample run. The two should be identical; otherwise the
message \texttt{Test failed!} will be displayed (on its own on a line) and
the target fails, so that \texttt{make} exits with an error and names the
test.

\vfill\pagebreak 

\begin{verbatim}
> make test1once
\end{verbatim}
runs the same curves again, but with the input fields of the argument
structure set once before the loop over the curves instead of once per
curve, as a program using the library may do since version~3.0.0
(Section~\ref{sec:library}); \texttt{rptest} prints a line whenever the
library has changed one of those fields, and the output is compared
against \texttt{testbase} as before.

\begin{verbatim}
> make test1many
\end{verbatim}
does what \texttt{make test1} does, but with a second collection of curves, in
\texttt{testdata-many.h}, whose reference output is \texttt{testbase-many}.
These curves have (comparatively) many rational points, which makes them
sieve quite differently from the random curves in \texttt{testdata.h}: the
smallest primes say little or nothing about them, so the first sieving stage
needs about sixteen moduli where a random curve needs about ten to let the
same proportion of numerators through (about one in three thousand), and
five times as many candidates per denominator survive the second stage; the
best division of labour between the stages is therefore a different one. The
two tests take about the same time and
between them cover the two regimes that occur in practice; both should be used
whenever the parameters are retuned, since a setting that suits one of them can
be a poor choice for the other. That is what \texttt{make tune} does with them;
see Section~\ref{sec:tuning}.

\begin{verbatim}
> make testhigh
> make testhighmany
\end{verbatim}
run the same two regimes at a height bound of $2\cdot10^5$ rather than the
$16383$ of \texttt{test1} and \texttt{test1many}, and take about a minute
each. The height bound decides a good deal of what a test measures. A sieve table
is built the first time a given prime meets a given class of denominators and
is reused afterwards, so the total cost of building the tables is bounded by
the number of primes and does not grow with the height, whereas the sieving
does: on the random curves the sieve tables are $4\%$ of \texttt{make test1}
but $0.3\%$ of \texttt{make testhigh}, the whole of the set-up per
denominator $8.5\%$ against $1.3\%$, and the two sieving stages themselves go
from $57\%$ to $90\%$ of the run. These are therefore the tests to judge a
change to the sieving loops by. \texttt{testhigh} uses the same curves as
\texttt{test1}; none of them has a rational point of height between $16383$
and $2\cdot10^5$, so the expected output is \texttt{testbase} again.
\texttt{testhighmany} uses the thirty
curves of \texttt{testdata-many.h} with the most rational points, on which
the small primes say least and the program has to look furthest up the table
of primes, collected in \texttt{testdata-high-many.h} with reference output
\texttt{testbase-high-many}.

\begin{verbatim}
> make test2
\end{verbatim}
calls \texttt{ratpoints} on a curve with lots of rational points and with
a fairly large height bound, times it, and compares the output to what is
expected.

\begin{verbatim}
> make test3
\end{verbatim}
runs the invocations of \texttt{ratpoints} collected in \texttt{test3.sh}
(which check that the bugs fixed in version~2.2.4, and those found while
version~3.0.0 was being developed, stay fixed) and
compares the output against \texttt{testbase3}.

\begin{verbatim}
> make test4
\end{verbatim}
runs a few hundred invocations of \texttt{ratpoints} collected in \texttt{test4.sh},
chosen so that nearly every line of the code runs at
least once (\texttt{make coverage} below says how nearly): degrees $1$ to~$9$ and up to the maximum of~$100$,
odd and even degrees, square and non-square
leading and constant coefficients, every reason to reverse the polynomial,
curves with no points modulo some prime or no admissible numerator modulo~$64$,
forbidden divisors of every kind, every option, restricted denominator
ranges and search intervals, height bounds from $1$ to~$2^{63}-1$ (those come
with a denominator range of one and a narrow search interval, so that the
run stays short), and every error message.  The output is compared against
\texttt{testbase4}; the whole test takes about a second.  Before each invocation
the script prints its arguments, so that a failing test can be found from
the difference between \texttt{test4.out} and \texttt{testbase4}.  The
reference was checked independently of the sieve: \texttt{verify-test4.py}
reads it and, for every invocation that prints points, searches every
coprime pair~$(a, b)$ in the range by brute force, in arbitrary precision;
it can be run again whenever the reference changes (it needs Python~3, and
PARI/GP for the 19~searches, on two curves, that are too large for it).

\begin{verbatim}
> make testapi
\end{verbatim}
runs \texttt{rpapi}, a program that tests what the library offers and the
command-line program cannot reach: the checks that
\texttt{find\_points\_work} makes on its arguments and the error codes it
returns, a callback that stops the search or declines a point, the flags,
and a sequence of searches on one initialized structure
(Section~\ref{sec:library}); the output is compared against
\texttt{testbase-api}.

\begin{verbatim}
> make test4configs
\end{verbatim}
runs \texttt{test4} again on the library built with each of the other
compile-time switches: the register widths 64, 128 (both variants) and
512 (emulated); \texttt{RATPOINTS\_CHUNK=1}, one register in the first
stage; \texttt{USE\_LONG\_IN\_PHASE\_2}, the second stage one word at a
time; \texttt{PRIME\_SIZE} 7 and~9, the smaller and the larger prime table;
and composite moduli off and unbounded.  Each is
built in a directory of its own, so that the build in the working
directory is left alone (\texttt{test4-configs.sh}); the points found
must not depend on any of these, so every build is compared against the
same \texttt{testbase4}. It takes about a minute.
The script exits with status~1 when an output
differs from the reference and with~2 when a build failed; so does
\texttt{coverage.sh} below.

\begin{verbatim}
> make coverage
\end{verbatim}
measures how much of the code the tests execute, with~\texttt{gcov}: it
builds the sources instrumented and unoptimised in~\texttt{build-coverage/},
runs the tests through that build (all of~\texttt{make test} but \texttt{test2}
and~\texttt{timing}, which add nothing but minutes) and prints,
for each source file, the lines executed
and the branches taken at least once (\texttt{coverage.sh}); the
annotated sources are left in that directory.
When this was written the tests executed $99.5\%$ of the lines
of~\texttt{find\_points.c}, $96\%$ of~\texttt{sift.c}, all of~\texttt{init.c}
and~\texttt{sturm.c} and $99\%$ of~\texttt{main.c}; what is left is guards
against states no caller can produce, the paths that a failure to create a
thread or to get memory would take, arms that other compile-time settings
reach (\texttt{test4configs} runs those builds, but this measures one), and
a few conditions the arithmetic makes impossible.

\begin{verbatim}
> make testthreads
\end{verbatim}
runs \texttt{rptest}, \texttt{rptest-many} and \texttt{rptest-degrees} with
two to five sieving threads and \texttt{test3} and \texttt{test4} with
three and seven, against the same references as before, since what the
program prints must not depend on the number of threads; then
\texttt{test4-threads.sh} builds the library three more times with the
blocks into which the loop over the denominators is cut
(Section~\ref{sec:threads}) fixed at 1, 3 and 17 steps, and runs
\texttt{test4} and \texttt{rptest} on each with three and six threads,
since the output must not depend on where the blocks end either.
About a minute.

\begin{verbatim}
> make tsan
\end{verbatim}
builds the library and the test programs in \texttt{build-tsan/} with
\texttt{gcc}'s ThreadSanitizer and runs the threaded tests on them
(\texttt{tsan.sh}, which has the compiler flags): every output must match
its reference, and the sanitizer must not report a data race. A few minutes.

\vfill\pagebreak 

\begin{verbatim}
> make timing
\end{verbatim}
times \texttt{ratpoints} on some curve with height bound $400\,000$.
This may be useful for comparisons. To time it with four threads, say
\begin{verbatim}
> make timing RPOPTS='-t 4'
\end{verbatim}

Finally,
\begin{verbatim}
> make test
\end{verbatim}
runs \texttt{test1}, \texttt{test1once}, \texttt{test1many},
\texttt{testdegrees}, \texttt{test2}, \texttt{test3}, \texttt{test4},
\texttt{testapi} and \texttt{timing},
and fails at the end if any of them failed, so that a script can
tell whether the build passed. It takes about a minute.
The two large-height tests are not part of it; they are meant for judging a
change to the sieve rather than for checking that the build works. Neither
are \texttt{testthreads} and \texttt{tsan}, which take longer than the rest
together.

\subsection{Install}

If you like, you can install the library, executable and header file
on your system.

\begin{verbatim}
> sudo make install
\end{verbatim}

The executable is copied to \texttt{/usr/local/bin/}, the
library to \texttt{/usr/local/lib/}, and the header file to
\texttt{/usr/local/include/}. You can change the \texttt{/usr/local} prefix
via
\begin{verbatim}
> make install INSTALL_DIR=...
\end{verbatim}

\subsection{Debugging}

In case you found a bug and would like to find out where it comes from,
or if you just want to see exactly what the program is doing, you can
do

\begin{verbatim}
> make debug
\end{verbatim}

in the working directory. This will build an executable \texttt{ratpoints-debug},
which, when run, will dump loads of output on the screen, so it is best
to send the output to a file, which you can then study at leisure.

\subsection{Cleaning up}

In order to get rid of the temporary files, do

\begin{verbatim}
> make clean
\end{verbatim}

To remove everything except the files from the archive, use

\begin{verbatim}
> make distclean
\end{verbatim}

\subsection{Requirements}

You need a C compiler (like \texttt{gcc}, which is the default specified
for the \texttt{CC} variable in \texttt{Makefile}; if necessary, it can be
changed there). Since version~3.0.0, the code assumes that \texttt{long} has
64~bits and refuses to compile otherwise; that is the case on every 64-bit
Unix-like system (but not under the Windows ABI, where \texttt{long} has
32~bits even on a 64-bit machine). Version~2.2.4 is the last one that
accommodates a 32-bit \texttt{long}.

In addition to the standard libraries, the program also requires the
\texttt{gmp} (GNU multi-precision) library~\cite{gmp}. Sieving on several
threads needs POSIX threads, which every Unix-like system has;
\texttt{make THREADS=0} builds without them (Section~\ref{sec:build}).

\subsection{List of files}

{\sloppy%
The archive contains the following files.
\begin{itemize}
  \item \texttt{Makefile}
  \item \texttt{ratpoints.h} --- the header file for programs using ratpoints
  \item \texttt{rp-private.h, primes.h} --- header files used internally
  \item \texttt{gen\_find\_points\_h.c, gen\_init\_sieve\_h.c} ---
        short programs that write additional header files depending
        on the system configuration
  \item \texttt{sift.c, init.c, sturm.c, find\_points.c} ---
        the source code for the ratpoints library
  \item \texttt{main.c, rptest.c, rpapi.c} --- the source code for the
        \texttt{ratpoints}, \texttt{rptest} and \texttt{rpapi} executables,
        respectively
  \item \texttt{testdata.h}, \texttt{testdata-many.h},
        \texttt{testdata-high-many.h}, \texttt{testdata-degrees.h} ---
        the curves used by the test runs
  \item \texttt{testbase}, \texttt{testbase-many},
        \texttt{testbase-high-many}, \texttt{testbase-degrees},
        \texttt{testbase2}, \texttt{testbase3}, \texttt{testbase4},
        \texttt{testbase-api} --- the output they are expected to produce
  \item \texttt{verify-test4.py} --- the brute-force check of
        \texttt{testbase4} (needs Python~3)
  \item \texttt{test3.sh}, \texttt{test4.sh} --- the scripts that
        \texttt{make test3} and \texttt{make test4} run; they have to be
        executable, as have the next five
  \item \texttt{test4-configs.sh}, \texttt{coverage.sh} --- what
        \texttt{make test4configs} and \texttt{make coverage} run
  \item \texttt{test4-threads.sh}, \texttt{tsan.sh} --- the second half of
        \texttt{make testthreads}, and what \texttt{make tsan} runs
  \item \texttt{tune.sh} --- what \texttt{make tune} runs to measure the
        machine-dependent constants (Section~\ref{sec:tuning})
  \item \texttt{bench\_init.c}, \texttt{bench\_check.c} --- benchmarks for
        two parts that are timed on their own: building the sieving tables,
        and the exact test.  Neither is part of the library; the second is
        how to remeasure the constants of Section~\ref{sec:sp3choice}
  \item \texttt{ratpoints-doc-3.1.tex} --- the source of this documentation
        (\texttt{make all} or \texttt{make doc} produce the PDF file)
  \item \texttt{README.md} --- a short description and an overview of what
        is new in this version
  \item \texttt{gpl-2.0.txt} --- the GNU license that applies to this program.
\end{itemize}%
}


\section{How to use \texttt{ratpoints}}

\subsection{Basic operation}

Let
\[ C : y^2 = a_n x^n + a_{n-1} x^{n-1} + \dots + a_1 x + a_0 \]
be your
curve, and let $H$ be the bound for the denominator and absolute value
of the numerator of the $x$-coordinate of the points you want to find.
The command to search for the points is then

\verb+> ratpoints '+$a_0$ $a_1$ \dots $a_{n-1}$ $a_n$\verb+' +$H$

The first argument to \texttt{ratpoints} is the list of the coefficients,
which have to be integers, are separated by spaces, and are listed
starting with the constant term. There is no bound on the size of the
coefficients; they are read as multi-precision integers.

The degree $n$ of the polynomial is limited to \texttt{RATPOINTS\_MAX\_DEGREE}
(defined in \texttt{ratpoints.h}), which is set to~$100$ by default.

The program will give an error message when the polynomial (considered
as a binary form of the smallest even degree $\ge n$) is not squarefree
(e.g., if you specify two leading zero coefficients).

The following subsections discuss how to modify the standard behavior.
This is achieved by adding options after the two required arguments.
Some of the options have arguments, others don't; the ordering of
the options and option-argument pairs is largely arbitrary.
The exceptions are options with opposite effect, where the last one given counts,
and the specification of the search intervals, which must be given in increasing order.

\subsection{Early abort}

By default, \texttt{ratpoints} will search the whole region that you specify
and print all the points it finds. If you just want to find {\em one}
point (for example when dealing with 2-covering curves of elliptic curves),
you can tell the program to quit after it has found one point via the
\verb+-1+ option, e.g.,

\verb+> ratpoints '-18 116 48 -12 30' 60000000 -1+

If you don't want to see points at infinity, specify the \verb+-i+ option.
The two options can be combined: \verb+-i -1+ will stop the program after
it has found one finite point.

\subsection{Changing the output}

By default, \texttt{ratpoints} prints some general information before and
after the points, which are given in the form
\verb+(+$x$\verb+ : +$y$\verb+ : +$z$\verb+)+, one per line.
Here, $(x : y : z)$ are the coordinates of the point considered as a point
in the $(1, \lceil n/2\rceil, 1)$-weighted projective plane, which is the
natural ambient space for the curve. The coordinates are integers with
$x$ and~$z$ coprime and $z > 0$, or $z = 0$ and $x = 1$ (for points at
infinity).

There are various ways to change this behavior.
\begin{itemize}\setlength{\parindent}{0mm}\setlength{\parskip}{1ex}
  \item To suppress all output except the points, use \verb+-q+ (for
        {\em quiet}).
  \item To add some more output explaining what the program is doing,
        use \verb+-v+ (for {\em verbose}). This has no effect if
        the \verb+-q+ option is also given.
  \item To suppress printing of the points, use \verb+-z+.
  \item If you only want to list the $x$-coordinates of point pairs rather
        than individual points, use~\verb+-y+.
  \item There are four options that influence the format in which the
        points are printed: \\
        \strut\quad \verb+-f+ {\em format} \verb+-fs+ {\em string-before}
        \verb+-fm+ {\em string-between} \verb+-fe+ {\em string-after} \\
        The arguments to all of them are strings; \verb+"\n"+, \verb+"\t"+,
        \verb+"\\"+ and \verb+"\%"+ are recognized and do what you expect.
        The {\em format} string can contain markers \verb+%x+, \verb+%y+,
        \verb+%z+ that will be replaced with the $x$, $y$, and $z$-coordinate
        of the point, respectively. The defaults are empty for
        {\em string-before}, {\em string-between} and {\em string-after},
        and \verb+"(%x : %y : %z)\n"+ for {\em format} when \verb+-y+
        is not specified; otherwise \verb+"(%x : %z)\n"+.

        The effect is as follows. Before any point is printed,
        {\em string-before} is output. Then every point is printed
        according to {\em format}. Between any two points, {\em string-between}
        is output, and after the last point, {\em string-after} is output.

        As an example, consider the effect of \\
        \strut\quad\verb+-f "[%x,%y,%z]" -fs "{" -fm "," -fe "}\n"+
\end{itemize}

\subsection{Restricting the search domain}

There are two ways to restrict the domain of the search. The first is
to restrict the range of denominators considered. This is done via \\
\strut\quad\verb+-dl+ $d_{\text{min}}$ \verb+-du+ $d_{\text{max}}$\,. \\
For example, to only look for integral points, you can say \verb+-du 1+.
By default, the lower limit is~$1$ and the upper limit is~$H$.

The other way is to restrict the search to a union of (closed) intervals.
If you only want to search for points in
\[ [l_1, u_1] \cup [l_2, u_2] \cup \dots \cup [l_k, u_k] \,, \]
you can specify this in the form \\
\strut\quad\verb+-l+ $l_1$ \verb+-u+ $u_1$ \verb+-l+ $l_2$ \verb+-u+ $u_2$ \dots
\verb+-l+ $l_k$ \verb+-u+ $u_k$ \,. \\
The first \verb+-l+ option is optional; $l_1$ defaults to~$-\infty$.
Similarly the last \verb+-u+ option is optional, and $u_k$ defaults
to~$+\infty$. The number $k$ of intervals is bounded by
\texttt{RATPOINTS\_MAX\_DEGREE} (usually, $100$).

\subsection{Using several threads}\label{sec:threads}

\strut\quad\verb+-t+ $n$ \\
runs the sieve on $n$ threads; \verb+-t 0+ uses as many as the machine has
processors. By default the program uses one thread, its own. With
\verb+-t+~$n$ for $n \ge 2$, $n$ further threads do the sieving and the exact
tests, while the program's own thread hands out the work and prints the
points. The loop over the denominators is cut into blocks of consecutive
denominators, whose length depends only on the search and not on the
number of threads; the threads take the blocks in order, and the points
of a block are printed only when every earlier block has been printed.
So the output is the same, line for line, however many threads are used,
including the reports of~\verb+-v+ and what \verb+-1+ finds: the sieving
parameters are chosen and revised from the same counts on the same
schedule, and the first point is the first point of the sequential
program, not the first one found in time. At most $33$ threads sieve at
once, so \verb+-t 0+ on a machine with more processors uses~$33$; a
build for more cores is described in Section~\ref{sec:build}.

Threads pay for long runs. On a 16-core AMD Ryzen~9 9950X3D, a search
with a height bound of $4\cdot10^5$ (the curve of \texttt{make timing})
takes $1.50$~seconds on one thread and $0.80$, $0.45$, $0.24$ and $0.16$
on 2, 4, 8 and~16 (measured by \textbf{Drew Sutherland}); the two long
test suites, 1008 random curves and 30 point-rich curves at height bound
$2\cdot10^5$, run $5.7$ and $7.9$ times faster on 16 threads than on one,
since every curve has a set-up that one thread does. On the author's
laptop, whose power budget lowers the clock as more cores work, the same
search takes $1.9$~seconds on one thread and $1.0$, $0.7$ and $0.4$ on 2,
4 and~10; the sieve's work is divided without loss (measured in processor
cycles, two threads spend what one does), so what the elapsed time shows
beyond that is the machine. At a small height bound the set-up of a
search is a large share of the run and there is little to gain; many
curves at a small height bound are better searched by independent
searches running in parallel (Section~\ref{sec:library}).

\subsection{Setting the parameters for the sieve} \label{fd}

It is possible to change the number of primes that are used in the various
stages of the algorithm. This is done by the following options. \\
\strut\quad\verb+-p+ $M$ \verb+-N+ $N$ \verb+-n+ $n$ \\
This sets the number of (odd) primes that are considered for the sieve
to~$M$, the number of primes that are actually used for the sieve to~$N$,
and the number of primes used in the first sieving stage to~$n$. Giving~$M$
also fixes it: without \verb+-p+ the program considers
\texttt{RATPOINTS\_DEFAULT\_NUM\_PRIMES} of them to begin with and looks
further only for a curve that needs it (Section~\ref{sec:spchoice}), whereas
with \verb+-p+ it looks at exactly~$M$. The program
will, if necessary, reduce $M$, $N$, and $n$ (in this order) to ensure
that $n \le N \le M \le \text{\texttt{RATPOINTS\_NUM\_PRIMES}}$. The latter
is $53$ by default (the number of odd primes $< 2^8$), and $30$ or~$96$ if
\texttt{PRIME\_SIZE} is set to~$7$ or~$9$, etc.

$M$ has a price in memory, because the sieving tables for a prime~$p$ take
$p$ rows of $p$ bits and so grow with~$p^2$: the primes below $128$ need about
5\,MB, those below $1024$ about 1.7\,GB (with 256-bit registers;
twice that with 512-bit registers).
Only the primes that $M$ actually asks for are reserved, so raising
\texttt{PRIME\_SIZE} in order to have larger primes available costs nothing as
long as $M$ is left alone.

If \verb+-n+ or \verb+-N+ is left out, the corresponding number is chosen
from the curve, as described in Section~\ref{sec:spchoice}. The two constants
that govern that choice can themselves be set: \\
\strut\quad\verb+-r+ $x$ \verb+-R+ $e$ \\
sets the target number of numerators surviving the first sieving stage per
64-bit word to~$x$ (a real number; this is what decides~$n$), and the number
of primes added to~$n$ to obtain~$N$ to~$e$. The mnemonic is that \verb+-r+ and
\verb+-R+ are to the automatic choice what \verb+-n+ and \verb+-N+ are to
the explicit one. Both are machine-dependent; the compiled-in defaults
(determined on the author's machine as of September~2026)
are \texttt{RATPOINTS\_SURVIVORS\_PER\_WORD}
and \texttt{RATPOINTS\_SP2\_EXTRA} in \texttt{ratpoints.h}.

Three further options control the same choice; they make it possible to
switch off or to measure certain aspects of the selection of the moduli: \\
\strut\quad\verb+-U+ $u$ \quad\verb+-C+ $c$ \quad\verb+-A+ $a$ \\
The first sets the length of run, in 64-bit words of numerators, at which a
second-stage prime just pays for setting itself up; \verb+-R+ is then the
number of extra primes for a run much longer than that, and a shorter run
gets proportionately fewer (Section~\ref{sec:spchoice}); \verb+-U 0+ makes
\verb+-R+ a flat offset. The second sets what building one row of a sieving
table costs, in units of what one first-stage prime costs per word, which
is what lets the program prefer a smaller prime to a larger one of the same
quality; \verb+-C 0+ ranks the primes by quality alone. The third says how
much of the choice the program should revise while it runs
(Section~\ref{sec:adapt}): $a = 0$ nothing, $a = 1$, the default, the
number of primes in the third stage, and $a = 2$ the number in the second
stage as well. It revises nothing it did not choose in the first place, so
\verb+-n+ or \verb+-N+ has the same effect as \verb+-A 0+.

A third sieving stage runs after those two, on the numerators that survive
them. It tests each of them against further primes one at a time, computing
the position in the table of admissible residues rather than reading a bit
out of a sieving table, so it needs no table and can use primes for which
building one would be out of the question. How many primes it uses is
chosen from the curve, and \\
\strut\quad\verb+-P+ $k$ \\
overrides that with exactly $k$ of them; \verb+-P 0+ switches the stage off,
and \\
\strut\quad\verb+-Q+ $q$ \\
sets instead the machine constant the choice turns on, namely what carrying
one more prime through the denominators costs as a fraction of one exact
test. It is to \verb+-P+ what \verb+-r+ is to \verb+-n+.
The choice weighs what a prime removes, which is a fraction of the survivors,
against what it costs, which is a little work per survivor and a little more
per denominator, so it depends on how many survivors a denominator has left
by then --- many for a curve with many rational points, few for a random one,
and more when the second stage has been given few primes, as it is on a
short run. On the random test curves the stage uses two or three primes at a
height bound of~$16\,383$ and mostly none at~$2\cdot10^5$, where the second
stage has done the work already; on the curves with many rational points it
uses up to ten. Section~\ref{sec:sp3choice} has the rule; the two machine constants in it
are \texttt{RATPOINTS\_SP3\_PER\_SURVIVOR} and
\texttt{RATPOINTS\_SP3\_PER\_DENOM} in~\texttt{ratpoints.h}.

Both of those constants are fractions of what one exact test costs, which
depends on the curve: the test evaluates a binary form of
the given degree at one pair of integers and takes an integer square root, so
it costs more at a high degree, and more again when the coefficients are
large enough to push the arithmetic into another limb. The program estimates
this cost from the degree, the largest coefficient and the height bound, and \\
\strut\quad\verb+-W+ $w$ \\
overrides the estimate with $w$ rdtsc cycles. An rdtsc cycle is one tick of
the processor's time-stamp counter, which the x86 instruction~\texttt{rdtsc}
reads; on current processors that counter advances at a fixed rate close to
the nominal clock rate, whatever frequency the core is running at, so it
measures time --- in units of a fraction of a nanosecond --- rather than the
core's own cycles. It is the unit in which the program measures itself:
\texttt{make bench\_check} reports it, and
\texttt{RATPOINTS\_CHECK\_REFERENCE} in~\texttt{ratpoints.h} is expressed
in it. \verb+-W 306+ assumes that every curve costs what a degree-6 curve
with small coefficients costs (on the author's machine).
Section~\ref{sec:sp3choice} has the formula, and \verb+-v+ prints
what it gives for the curve in~hand.

There are two more parameters that can be set. \\
\strut\quad\verb+-F+ $D$ \\
sets the maximal number of `forbidden divisors'. If the degree is even,
the denominator of any $x$-coordinate of a rational point is restricted at
certain primes~$p$: it cannot be divisible by~$p$ when the leading
coefficient is a non-square mod~$p$, and when $p$ divides the leading
coefficient its valuation at~$p$ is restricted instead --- the program
works out which valuations can occur at all (Section~\ref{sec:elim}), which
typically rules out a power of~$p$, or one particular valuation, or~$p$
itself. Each such prime counts as one forbidden divisor. The program
constructs a list of them against which to check the denominators, taking
the primes up to the square root of the height bound, and this option caps
the length of the list; the default is~$64$, which the list reaches at a
height bound of about half a million. A prime in the list costs a few
instructions per 64 denominators. With \verb+-v+ the program prints what it found.

The other option is \\
\strut\quad\verb+-S+ $S$ \quad or\quad \verb+-s+ \,. \\
In the first form, it sets the number of refinement (interval halving)
steps in the isolation of the real components to~$S$. This part of the
program computes a Sturm sequence for the polynomial and uses it in order
to find a union of intervals that contains the intervals of positivity
of the polynomial. This is done by a successive subdivision of the
interval $]{-\infty}, +\infty[$. The number~$S$ sets the recursion depth.
$S$ can be omitted; then it is given a default value. The option
\verb+-s+ skips the Sturm sequence computation completely. This also has
the effect of removing most of the check for squarefreeness.

\subsection{Switching off optimizations}

The options \verb+-k+ and \verb+-j+ can be used to prevent the program
from reversing the polynomial, which it usually does if this will lead
to faster operation (\verb+-k+), and to prevent the program from using
the Jacobi symbol test on the denominators (\verb+-j+). This last test
extends the `forbidden divisors' method described in
Subsection~\ref{fd} above by computing the Jacobi symbol $(l/d)$, where
$l$ is the leading coefficient and $d$ is the odd and coprime-to-$l$
part of the denominator. If the symbol is~$-1$, the denominator need not
be considered. The use of these options may be questionable, unless you
want to see how much performance is gained by using these optimizations.

\subsection{Switching off exact testing of points}

The option \verb+-x+ will prevent the program from checking whether the
potential points that survive the sieve really give rise to rational
points, so the $x$-coordinates that are output may not actually correspond
to rational points. It also effectively sets the \verb+-y+ option, since the
$y$-coordinates are not computed.

\subsection{Overriding previous options}

The options \verb+-I+, \verb+-Y+, \verb+-Z+, \verb+-S+, \verb+-K+, \verb+-J+,
\verb+-X+ can be used to cancel the effect of the corresponding lower-case
option (and thereby restore the default behavior),
when this occurs earlier in the list of options. This may be useful if
you want to set a default behavior that is different from what the
program does out-of-the-box (e.g., in a shell script), but want the caller
to be able to override this change.


\section{How to use the library}\label{sec:library}

It is possible to use the ratpoints machinery from within your own programs.

The library \texttt{libratpoints.a} provides the following functions.

\begin{verbatim}
long find_points(ratpoints_args*,
                 int proc(long, long, const mpz_t, void*, int*),
                 void*);

void find_points_init(ratpoints_args*);

long find_points_work(ratpoints_args*,
                      int proc(long, long, const mpz_t, void*, int*),
                      void*);

void find_points_clear(ratpoints_args*);
\end{verbatim}

The passing of arguments to these functions is via the \texttt{ratpoints\_args}
structure, which is defined as follows.

\begin{verbatim}
typedef struct { mpz_t *cof; long degree; long height;
                 ratpoints_interval *domain; long num_inter;
                 long b_low; long b_high; long sp1; long sp2; long sp3;
                 double survivors_per_word; long sp2_extra; double sp2_u0;
                 double cost_table; long adapt;
                 long sp3_extra; double sp3_per_denom; double check_cost;
                 long array_size;
                 long sturm; long num_primes; long max_forbidden;
                 long num_threads;
                 unsigned int flags;
                 long sp1_used; long sp2_used; long sp3_used; ...}
        ratpoints_args;
\end{verbatim}

The dots at the end stand for additional fields that are used internally
and are not of interest here. All fields up to and including \texttt{flags}
are inputs, and they come back from a call as they went in, with three
exceptions that are described below: the coefficient array and
\texttt{degree}, which on return describe the polynomial the program
worked with; \texttt{num\_inter} and \texttt{domain}, which on return
describe the region it searched; and the flag bits that report on the
run. (Up to version 2.2.4 the program also wrote the values it chose or
normalised into the other input fields, so that a caller who fills the
structure once and then loops over curves had to reset them before every
call; this is no longer necessary.) Every input field has to be set before
the call; a negative value asks for the compiled-in default, or for the
value to be chosen from the curve, as described below.
The three \texttt{\_used} fields are outputs; see below.

When calling \texttt{find\_points}, the
\texttt{cof} field must point to an array of size $\text{\texttt{degree}} + 1$
of properly initialized gmp integers; these are the coefficients of
the polynomial. The program can alter the values of the integers in this
array: it may reverse the polynomial, and it drops a zero constant term
on that occasion, in which case \texttt{degree} is lowered by one so that
\texttt{cof} and \texttt{degree} still describe the polynomial the program
worked with (the \texttt{RATPOINTS\_REVERSED} flag says when this has
happened). Leading zero coefficients are dropped in the same way, lowering
\texttt{degree}, without a flag.
To prevent the reversal, set the \texttt{RATPOINTS\_NO\_REVERSE} bit
in \texttt{flags}; this may result in a loss of performance, though.

\texttt{height} gives the height bound. It is an error for the
degree or the height bound to be nonpositive; in this case the function
returns the value \texttt{RATPOINTS\_BAD\_ARGS}.

The field \texttt{domain} must contain a pointer to an array of
\texttt{ratpoints\_interval} structures of length at least \texttt{num\_inter}
plus \texttt{degree}. This array gives (in its first \texttt{num\_inter} entries)
the intervals for the search region. The type \texttt{ratpoints\_interval}
is just

\begin{verbatim}
typedef struct {double low; double up;} ratpoints_interval;
\end{verbatim}

the meaning of this should be clear. In \texttt{ratpoints}, this is set
via the \verb+-l+ and \verb+-u+ options. Usually, you do not want to
restrict the range of $x$-coordinates; then you set \texttt{num\_inter = 0}
(but you still have to fill \texttt{domain} with a valid pointer to at least
\texttt{degree} intervals; a null pointer makes the call return
\texttt{RATPOINTS\_BAD\_ARGS}, since the program stores the default
interval $]-\texttt{height}, \texttt{height}[$ there. If \texttt{sturm} is
negative, one interval is enough.)
The program intersects the given intervals with the region where $f$ is
positive (unless \texttt{sturm} has a negative value, which may lead to a
performance loss) and searches the result; on return, \texttt{num\_inter}
and the first \texttt{num\_inter} entries of the array describe that
result, which is why the array has to be longer than the input. So
\texttt{num\_inter} (and \texttt{domain}, if you restrict the search) must
be set before every call.

\texttt{b\_low} and \texttt{b\_high}
carry the lower and upper bounds for the denominator; if non-positive, they are
set to~$1$ and the height bound, respectively. In \texttt{ratpoints}, these
fields are set by the \verb+-dl+ and \verb+-du+ options.

\texttt{sp1} and \texttt{sp2}
specify the number of primes to be used in the first sieving stage and
in the two first sieving stages together, respectively. If negative, they are chosen
from the curve, as described in Section~\ref{sec:spchoice}; this is what the
command line program does when \verb+-n+ and \verb+-N+ are absent. In
\texttt{ratpoints}, these fields are set by the
\verb+-n+ and \verb+-N+ options. What was actually used, chosen or not,
is reported in \texttt{sp1\_used} and \texttt{sp2\_used}. The fields
\texttt{survivors\_per\_word} and \texttt{sp2\_extra} hold the two
constants that govern that choice (options \verb+-r+ and~\verb+-R+); if
negative, the values compiled into \texttt{ratpoints.h} are used. So do
\texttt{sp2\_u0} and~\texttt{cost\_table} (options \verb+-U+ and~\verb+-C+),
which say how long a run has to be for a second-stage prime to pay for
setting itself up and what one row of a sieving table costs; zero in either
switches off the corresponding refinement. \texttt{adapt} set to zero stops
the program revising these numbers as it runs (option~\verb+-A+,
Section~\ref{sec:adapt}); it revises them only when it chose them, so a
non-negative \texttt{sp1} or~\texttt{sp2} has the same effect.

\texttt{sp3\_used}
reports the number of primes used in all three stages together; the third
stage used the difference between that and \texttt{sp2\_used}.
(\texttt{sp3} itself is a working field of no interest to the caller.)
What goes in is \texttt{sp3\_extra}, which fixes that difference when
it is non-negative (option~\verb+-P+) and otherwise leaves the number to be
chosen from the curve, and~\texttt{sp3\_per\_denom} (option~\verb+-Q+),
which is the machine constant that choice turns on and takes the value
compiled into \texttt{ratpoints.h} when negative. \texttt{check\_cost}
(option~\verb+-W+) says what one exact test costs, in rdtsc cycles
(Section~\ref{fd}), and both
of those constants are fractions of it; when negative, which is the usual
case, the program estimates it from the degree, the coefficients and the
height bound as described in Section~\ref{sec:sp3choice}.
Similar statements are true for~\texttt{num\_primes} (option~\verb+-p+),
\texttt{sturm} (options \verb+-s+, \verb+-S+) and~\texttt{max\_forbidden}
(option~\verb+-F+). Note that a non-negative \texttt{num\_primes} is a hard
limit: only when it is negative does the program look at further primes of
its own accord for a curve that needs them. The field~\texttt{array\_size} specifies the maximal size
(in bit arrays; a bit array
contains 64, 128, 256 or~512 bits, depending on the configuration) of the
array that is used in the first sieving stage. If non-positive, it is
set to a default value. The various default values are defined at the
beginning of the header file~\texttt{ratpoints.h}, where also the maximal
degree of the polynomial is set.

\texttt{num\_threads} says how many threads sieve (option~\verb+-t+,
Section~\ref{sec:threads}): $1$ or~$0$, the calling thread alone; $n \ge 2$,
$n$ threads besides the calling one, which hands out the work and delivers
the points; negative, as many threads as there are online processors. The
threads are created at the first call that asks for them and kept, with
their working memory, until \texttt{find\_points\_clear}; so a program
that searches many curves should use the \texttt{find\_points\_init},
\texttt{find\_points\_work}, \texttt{find\_points\_clear} sequence
described below rather than \texttt{find\_points}, which creates and ends
them on every call. A call with another value of \texttt{num\_threads}
replaces them. If the threads cannot be created, the calling thread sieves
alone. No more than $33$ threads are created, whatever the field says,
unless the library is built for more (Section~\ref{sec:build}). A library
built without threads (\texttt{make THREADS=0}) ignores the field.

The \texttt{flags} field holds a number of bit flags.
\begin{itemize}
  \item \texttt{RATPOINTS\_NO\_CHECK} --- when set, do not check whether
        the surviving $x$-coordinates give rise to rational points
        (set by the \verb+-x+~option to~\texttt{ratpoints}).
  \item \texttt{RATPOINTS\_NO\_Y} --- only list $x$-coordinates (in the
        form $(x : z) \in \BP^1(\Q)$) instead of actual points (with
        a $y$-coordinate); this is set by the \verb+-y+~option
        to~\texttt{ratpoints}.
  \item \texttt{RATPOINTS\_NO\_REVERSE} --- when set, do not allow reversal
        of the polynomial (set by the \verb+-k+~option to~\texttt{ratpoints}).
  \item \texttt{RATPOINTS\_NO\_JACOBI} --- when set, prevent the use of the
        Jacobi symbol test (set by the \verb+-j+~option to~\texttt{ratpoints}).
  \item \texttt{RATPOINTS\_VERBOSE} --- when set, causes the procedure to
        print some output on what it is doing (set by the \verb+-v+~option
        to~\texttt{ratpoints}).
\end{itemize}

There are some other flags that are used internally. One of them might be
of interest:
\begin{itemize}
  \item \texttt{RATPOINTS\_REVERSED} --- when set after the function call,
        this indicates that the polynomial has been reversed (and the
        contents of the \texttt{cof} array have been modified).
\end{itemize}

The main vehicle for passing information back to the caller is the
\texttt{proc} function argument together with the pointer~\texttt{info}.
This function \\
\strut\quad\verb+int proc(long x, long z, const mpz_t y, void *info, int *quit)+ \\
is called whenever a point was found. \texttt{x}, \texttt{y} and~\texttt{z} are
the coordinates of the point (where \texttt{y} is a gmp integer).
\texttt{info} is the pointer that was passed to~\texttt{find\_points}; this
can be used to store information that should persist between calls
to \texttt{proc}. If \texttt{*quit} is set to a non-zero value, this indicates
that \texttt{find\_points} should abort the point search and return immediately;
otherwise the search continues. This is how the \verb+-1+~option works.
The return value is taken as a weight
for counting the points; usually it will be~1.

Whatever \texttt{num\_threads} says, \texttt{proc} is called from the
thread that called \texttt{find\_points}, never from a sieving thread, and
for the same points in the same order as a single thread finds them; the
return value of \texttt{find\_points} and the point at which \texttt{*quit}
takes effect are the same as well. So \texttt{proc} and \texttt{info} need
no protection against concurrent calls. A sieving thread keeps the points
it finds until every earlier block of the search has been delivered; if it
cannot get the memory for that, the search fails with the return value
\texttt{RATPOINTS\_NO\_MEMORY} rather than deliver fewer points than it
found.

The usual framework for using \texttt{find\_points} is as follows.

\begin{verbatim}
(...)
#include "ratpoints.h"

(...)

mpz_t c[RATPOINTS_MAX_DEGREE+1];  /* The coefficients of f */
ratpoints_interval domain[2*RATPOINTS_MAX_DEGREE];
                     /* This contains the intervals representing the
                        search region */

/*****************************************************************
 * function that processes the points                            *
 *****************************************************************/

typedef struct {...} data;

int process(long x, long z, const mpz_t y, void *info0, int *quit)
{ data *info = (data *)info0;

  (...)

  return(1);
}
\end{verbatim}

\begin{verbatim}
/*****************************************************************
 * main                                                          *
 *****************************************************************/

int main(int argc, char *argv[])
{
  long total, n;
  ratpoints_args args;

  long degree        = 6;
  long height        = 16383;
  long sieve_primes1 = -1; /* negative: chosen from the curve */
  long sieve_primes2 = -1; /* negative: chosen from the curve */
  long num_primes    = -1; /* negative: the default, more if needed */
  long max_forbidden = -1; /* negative: the compiled-in default */
  long num_threads   = 1;  /* sieve on this thread; see -t */
  long b_low         = 1;
  long b_high        = height;
  long sturm_iter    = RATPOINTS_DEFAULT_STURM;
  long array_size    = RATPOINTS_ARRAY_SIZE;
  int no_check       = 0;
  int no_y           = 0;
  int no_reverse     = 0;
  int no_jacobi      = 0;
  int no_output      = 0;

  unsigned int flags = 0;

  data *info = malloc(sizeof(data));

  /* initialize multi-precision integer variables */
  for(n = 0; n <= degree; n++) { mpz_init(c[n]); }

  (...)

  { /* set up polynomial */
    long k;
    for(k = 0; k < 7; k++) { mpz_set_si(c[k], ...); }

    args.cof           = &c[0];
    args.degree        = 6;
    args.height        = height;
    args.domain        = &domain[0];
    args.num_inter     = 0;
    args.b_low         = b_low;
    args.b_high        = b_high;
    args.sp1           = sieve_primes1;
    args.sp2           = sieve_primes2;
    args.survivors_per_word = -1.0; /* negative: the compiled-in default */
    args.sp2_extra     = -1;
    args.sp2_u0        = -1.0;
    args.cost_table    = -1.0;
    args.adapt         = -1;
    args.sp3_extra     = -1;        /* negative: chosen from the curve */
    args.sp3_per_denom = -1.0;
    args.check_cost    = -1.0;
    args.array_size    = array_size;
    args.sturm         = sturm_iter;
    args.num_primes    = num_primes;
    args.max_forbidden = max_forbidden;
    args.num_threads   = num_threads;
    args.flags         = flags;

    info->... = ...;
    (...)

    total = find_points(&args, process, (void *)info);
    if(total == RATPOINTS_NON_SQUAREFREE)
    { ... }
    if(total == RATPOINTS_BAD_ARGS)
    { ... }

    (...)

  }

  /* clean up multi-precision integer variables */
  for(n = 0; n <= degree; n++) {mpz_clear(c[n]); }

  return(0);
}
\end{verbatim}

If points are to be searched on many curves of the same degree,
then it is slightly more efficient to use the sequence

\begin{verbatim}
  args.degree = degree; /* this information is needed */
  find_points_init(&args);

  for( ... )
  { ...
    total = find_points_work(&args, process, (void *)info);
    ...
  }

  find_points_clear(&args);
\end{verbatim}

This avoids the repeated allocation and freeing of memory. Inside the loop
only \texttt{cof}, \texttt{degree} and~\texttt{num\_inter} (and whatever
else changes from curve to curve) need to be set; the other input fields
keep their values across calls. \texttt{rptest.c} with
the \verb+-O+~option does exactly this, and \texttt{make test1once} checks
that the points found are the same as with the fields set once per curve.

Searches on different \texttt{ratpoints\_args} structures are independent
of each other: the library keeps no state outside the structure (its
tables of primes and residues are constant), so a program may run several
searches at the same time on threads of its own, each with its own
structure, callback and \texttt{info}, and each with a
\texttt{num\_threads} of its own. For many searches at a small height
bound this is the better way to use a machine's cores
(Section~\ref{sec:threads}).

For practical examples, see \texttt{main.c} (the code that wraps
\texttt{find\_points} for the command line program~\texttt{ratpoints}) or
\texttt{rptest.c} (which runs \texttt{find\_points} on some test data).


\section{Implementation} \label{impl}

\subsection{Overview}

Let $F(x, z)$ be the binary form of even degree corresponding to the
polynomial on the right hand side of the curve equation. The basic idea
is to let $b$ run from~$1$ to the height bound~$H$, and for each~$b$, let
$a$ run from $-H$ to~$H$, and finally, for each coprime pair~$(a,b)$ check
if $F(a, b)$ is a square.

Of course, in this form, this would take a very long time. To speed up the
process, we try to eliminate quickly as many pairs $(a, b)$ as possible
before the actual test. This can be done by `quadratic sieving' modulo
several primes: if $F(a, b)$ is a square, it certainly has to be a square
mod~$p$, and so we can rule out all $(a, b)$ that do not satisfy this
condition. Another ingredient is to represent (for a fixed~$b$) the
various values of~$a$ by bits and treat all the bits in a `bit array'
(of 64 bits, or of 128, 256 or 512 bits when SSE/AVX instructions are
used) in parallel. For this, we organize the sieving information for each
prime~$p$ into $p$~arrays of $p$~words each, one such array for every $b$
mod~$p$, such that the $j$th~bit in this array is set if and only if
$F(j, b)$ is a square mod~$p$.

For each `denominator' $b$, we then set up an array of words whose bits
represent the range \hbox{$-H \le a \le H$} (aligned so that $a = 0$ corresponds
to the $0$th bit of a word); the bits are initially set. Then for each
of the sieving primes~$p$, we perform a bit-wise {\em and} operation
between this array and the sieving information at~$p$. In a first stage,
this is done on the whole array; after this first stage, each remaining
(`surviving') non-zero word in the array is subject to tests with more
primes. If some bits are still set after this second stage of sieving,
the corresponding pairs $(a, b)$ (if coprime) are tested against a few
further primes one at a time, without sieving tables (the third stage),
and those that pass are then checked exactly.

In the following subsections, we discuss a number of improvements that
were made.

\subsection{Sorting the primes}

This idea is due to \textbf{Colin Stahlke}. We do not just take the first so many
primes in increasing order for the sieving, but we first compute the
number of points the curve has mod~$p$ for a number of primes~$p$ and then
sort the primes according to the fraction of $x$-coordinates that give points.
We then take those primes for the sieving that have the smallest fraction
of `surviving' $x$-coordinates. In this way, we need fewer sieving primes
to achieve a comparable reduction of point tests. Since version~3.0.0 the
candidates include the odd composite numbers below~$64$
(a candidate sharing a prime with one already taken is passed over),
and the ranking also takes into account what they cost to set up
(Section~\ref{sec:spchoice}).

\subsection{Using connected components}

We note that we can only have points when $F(a, b)$ is non-negative.
If we can determine intervals on which $f(x) = F(x, 1)$ is negative, then
we do not have to look for points in these intervals. The necessary
computations can be performed exactly, by computing a Sturm sequence
for~$f$ and counting the number of sign changes at various points,
see~\cite[Thm.~4.1.10]{Cohen}. This can in particular tell us whether
$F$ is negative definite, in which case we have already proved that there
are no rational points in the curve. If there are real zeros, we use
a subdivision method in order to find a collection of intervals containing
the projections of the connected components of the curve over~$\R$.

\subsection{Using 2-adic information}

\textbf{John Cremona} suggested that in some cases, one can determine beforehand
that all points will have odd numerators~$a$, and so we can pack the
bits more tightly by only representing odd numbers. This approach can be
extended. We first find all solutions mod~$64$. Then for every
residue class mod~$64$ of the denominator~$b$, we can find the residue
classes mod~$64$ of potential numerators~$a$. If there are none, then we
can eliminate $b$ as a denominator altogether. If all potential $a$'s lie
in one residue class modulo~$2^k$, we can restrict the sieving to such numerators
($2^k = 2$, $4$ and $8$ occur on the test curves), and the arrays
are shorter by a factor of~$2^k$.

\subsection{Elimination of denominators}

Depending on the equation, certain denominators can be excluded. We have
seen an instance of this in the previous subsection, but we can also work
with odd primes.

\subsubsection{Odd degree and $\pm$monic}

If the polynomial has odd degree and leading coefficient $\pm 1$, then
the denominator has to be a square. This reduces the time complexity
tremendously, from $O(H^2)$ to~$O(H^{3/2})$.

\subsubsection{Odd degree general}

In general, when $f$ has odd degree, the denominator has to be `almost'
a square: it must be a square times a (squarefree) divisor of the leading
coefficient, and there are further restrictions on the parity of the
valuation at~$p$, for $p$ dividing the leading coefficient, when this
valuation is sufficiently large. This gives the same type of time complexity
as in the monic case, but with a larger constant.

\subsubsection{Even degree} \label{sec:elim}

Let $l$ be the leading coefficient. If $l$ is not a square, we can exclude
denominators divisible by a (small) prime~$p$ such that $l$ is a non-square
mod~$p$: for $p \mid b$ we have $F(a,b) \equiv l a^n \bmod p$ with $p \nmid
a$, which is a non-square as well. In addition, we use the necessary
condition that the Jacobi symbol $\left(\frac{l}{b'}\right)$ must be $+1$,
where $b'$ is the part of $b$ coprime to~$2l$. The Jacobi symbol
is computed (for most denominators, at least) as a product
of Legendre symbols over the prime divisors of~$l$; they are determined
by a fast table lookup. When $l$ has a prime factor above~$1021$, or more odd
prime factors than the tables allow for (sixteen, or 8~KB), or the
denominators exceed~$2^{32}$, the computation falls back on a binary algorithm.

When $p$ divides~$l$, we can instead consider $p$-adic valuations.
Write $v_p(b) = m \ge 1$ and let $w_j$ be
the valuation at~$p$ of the coefficient of~$b^j$ in~$F$, so that the term
with~$b^j$ has valuation $w_j + jm$. If one term has strictly smaller
valuation than all the others, it determines $v_p(F(a,b))$ and the residue
of $F(a,b)/p^{v}$ mod~$p$. If that valuation is odd, no denominator with
$v_p(b) = m$ can occur. If it is even and $j$ is even, then $a^{n-j}$ and
$(b/p^m)^j$ are both squares mod~$p$, so none can occur when the unit part
of that coefficient is a non-square mod~$p$.
For large~$m$ the leading term wins by itself, so when $v_p(l)$ is
odd, or even with a non-square unit part, every sufficiently large valuation
is excluded: the denominator cannot be divisible by a certain power of~$p$.
For example, when $m = 1$, then $p$ itself is excluded as a divisor
when $p$ divides the next coefficient, and $p^2$ is excluded otherwise.

While computing the sieving information, we make a list of these primes,
each with the set of valuations it excludes, which we then use to eliminate
denominators: a prime that is excluded altogether is tested for with bit
arrays in the manner of the sieve, the others by computing the valuation of
the denominator. The bit arrays cover every excluded prime up to the square
root of the upper bound for the denominators (but at most up to~$1021$, the
largest prime in the program's list, and at most as many primes as
\verb+-F+ allows), building the patterns for the
primes beyond the compiled sieving primes as needed.
Note that past the square root a denominator can have at most one excluded
prime left, and the Jacobi symbol sees that one.

\subsection{Reversing the polynomial}

The exclusion of denominators described above is more or less effective,
depending on the situation; it is the better the earlier the case was
described. Since the height of the points does not change under the
transformation $(x,y) \to (1/x, y/x^{\lceil n/2 \rceil})$, we can as well search on the
reversed polynomial $F(z,x)$. This will result in a speed-up if the
reversed polynomial belongs to a `better' class than the original.

\subsection{Some general remarks}\label{sec:remarks}

Much of what was learnt about where the time goes on a current processor
while this version was developed ran against expectation.

Integer division is slow,
and the earlier versions avoided it wherever they reasonably
could. However, division turned out not to be where the time goes --- the
integer divider is busy for less than $0.6\%$ of the cycles --- and the
subtraction chains written to avoid it were not that fast either: their
branches are predicted well, but five dependent conditional subtractions
take some twenty instructions more than two multiplications and a shift.
They have been replaced by reductions using a multiplication with a
precomputed reciprocal, which gained two to three per cent.

What does cost is memory traffic and control flow. The innermost loop of
the first sieving stage loads a row of a sieving table for each modulus and
\texttt{AND}s it into the registers, and it runs at the rate at which the
first-level cache delivers those rows: an \texttt{AND} costs about one cycle
as long as the rows of all first-stage moduli together fit into that cache
(some 25\,KB), and more for every kilobyte beyond. This is why the composite
sieving moduli are capped at~$64$ although larger ones would remove more
numerators per \texttt{AND}, why every row is stored \texttt{RBA\_PACK}
times over so that each load is aligned, and why nothing should be added to
the structure this loop reads (see Appendix~\ref{app:notes}).

A mispredicted branch costs fifteen to twenty cycles, the time of as many
\texttt{AND}s, and a branch that depends on the data mispredicts often.
This effect was particularly pronounced in previous versions at small height bounds.
The relevant parts (e.g., the Jacobi symbol, the gcd, locating the surviving
bits) now avoid data-dependent branches, resulting in a reduction of
the number of mispredicted branches to between $3\%$ and~$40\%$ of
what version~2.2.4 had.

For more details and further lessons, see Appendix~\ref{app:notes}.


\section{Fine-tuning the parameters}\label{sec:tuning}

Most of what the sieve does is decided from the curve before a single
numerator is sieved: which moduli to sieve with and how many in each of the
three stages, how long the run is going to be, and what one exact test will
cost. The rules weigh running-time costs against each other, so they rest on
a few constants that describe the machine rather than the curve. Reasonable
values (based on tests performed on the author's machine)
are compiled into \texttt{ratpoints.h}. \texttt{make tune} measures
better ones (Section~\ref{sec:retune}), and each of them can be set on the
command line for a single run. This section says what the rules are and
what the constants mean; Section~\ref{impl} says how the sieve works, and
Appendix~\ref{app:notes} records how the rules were measured and where each
constant came from.

\subsection{How the moduli of the first two stages are chosen}\label{sec:spchoice}

While collecting the sieving information, the program computes for each
modulus~$m$ --- the odd composite numbers below~$64$ and the primes in the
range considered --- the
density~$r_m$, the proportion of residues~$a$ modulo~$m$ for which $f(a)$
is a square modulo~$m$. The expected proportion of numerators surviving a set of
(pairwise coprime)
moduli is the product of their densities, so a modulus with a small density
gives more information.

The moduli are not ranked by density alone, because they do not cost the
same. A modulus~$m$ needs a sieving table of $m$~rows, built once for each
class of denominators modulo~$m$ that turns up --- at most $m$ times, however
long the run --- and an entry per denominator in the list of reduced
denominators; per word of numerators both fall like~$1/U$, where $U$ is the
number of words the run is going to sweep. The program computes~$U$ and the
number of denominators from the search region before it looks at a modulus,
and ranks the moduli by
\[ \frac{\text{cost}(m)}{-\log r_m} \,, \]
the cost per word divided by the information. At equal quality a smaller
modulus is preferred, and at a small height bound strongly so.
\texttt{RATPOINTS\_COST\_TABLE} (option~\verb+-C+) is what building one
row of a table costs, in units of one first-stage \texttt{AND} per word;
\verb+-C 0+ ranks by density alone.

The first stage is filled from the top of that ranking for as long as a
modulus pays: it is taken while the survivors per word it removes, weighted
by what a survivor costs in the stages after it, exceed
\texttt{RATPOINTS\_SURVIVORS\_PER\_WORD} (option~\verb+-r+) times what the
modulus costs per word --- one \texttt{AND}, plus its tables, its
per-denominator entries and the fetch of its row, spread over the run. That
threshold is the ratio of two machine costs and the most important constant
of all. With it the first stage stops at a few thousandths of a surviving
bit per word whatever the curve. The second stage then gets
\[ \texttt{sp2} = \texttt{sp1}
   + \frac{\texttt{RATPOINTS\_SP2\_EXTRA}}{1 + \texttt{RATPOINTS\_SP2\_U0}/U} \]
further moduli, capped by the number that carry any information at all.
\texttt{RATPOINTS\_SP2\_EXTRA} (\verb+-R+) is the number an arbitrarily long
run wants, and \texttt{RATPOINTS\_SP2\_U0} (\verb+-U+) is the length of run
at which setting a prime up costs as much as sieving with it, so that a
short run gets proportionately fewer; \verb+-U 0+ makes the offset flat. The
second stage has a ranking of its own, since a modulus
there is applied only to the bit arrays that survived the first stage, and
its table therefore weighs some three times more heavily against it; this is
why the second stage nearly always takes primes.

On a curve with very many rational points, $f(a)$ tends to be a square
modulo~$p$ for every residue~$a$ when $p$ is one of the smallest primes.
Such a prime says next to nothing: only that numerator and denominator are
not both divisible by~$p$, which is a density of $1 - 1/p^2$, and only when
denominators divisible by~$p$ occur at all (otherwise it is discarded). The
program keeps it as a candidate, since this little comes free when the prime
is a factor of a composite modulus --- $33$ in place of~$11$, say --- but does
not count it among the primes the rule wants,
and the thirty primes the program looks at by default may not leave enough
for the rule to have what it wants; the program then looks at further
primes, one at a time, until they suffice or the table ends. Every prime
added can only lower \texttt{sp1} and raise the number of informative
primes, so this stops of its own accord; random curves do not trigger it.
Giving \verb+-p+ explicitly switches it off: an explicit number of
primes is a hard limit, as are \verb+-n+ and \verb+-N+ for the first two stages.

\subsection{How many primes the third stage uses}\label{sec:sp3choice}

The third stage tests the survivors of the first two stages one at a time, against
further primes and without a sieving table. A prime saves the fraction
$1 - r_p$ of the survivors that reach it, each of which would otherwise go
to the exact test with its multi-precision arithmetic and its integer square
root; it costs a multiplication and a reduction for every survivor tested,
and one step per denominator to keep the denominator reduced modulo it,
whether or not a survivor turns up. So a prime is worth adding while
\[ S\,(1 - q_s - r_p) > q_d \,, \]
where $S$ is the expected number of survivors a denominator still has and
$q_s$, $q_d$ are the two costs as fractions of one exact test ---
\texttt{RATPOINTS\_SP3\_PER\_SURVIVOR} and \texttt{RATPOINTS\_SP3\_PER\_DENOM}
in \texttt{ratpoints.h}; option~\verb+-Q+ sets the second. Since $S$ falls
by the factor~$r_p$ with each prime, this stops of its own accord. On a
random curve the stage uses two or three primes on a short run, where the
second stage has been given few, and mostly none on a long one; on a curve
with many rational points it uses up to ten, and it
looks further up the table of primes whenever none in hand pays but a prime
of the density the curve has been offering would. \verb+-P+~$k$ fixes the
number at~$k$ instead, and \verb+-P 0+ switches the stage off.

$S$ is estimated from the number of numerators a denominator has to
consider, thinned by the 2-adic information, by the moduli of the first two
stages, and by the proportion of numerators coprime to the denominator ---
that last because the test for common factors runs before this stage.
This removes copies of $x$-coordinates that have already been dealt with,
in particular those coming from rational points, which would survive
every test.

What one exact test costs depends on the curve, not only on the machine. The
test evaluates the binary form $F(a,b)$ of degree~$d$ by Horner's rule and
takes an integer square root, so it costs a step per degree on numbers that
grow to about $\log_2|c| + d\log_2 H$ bits, plus a fixed amount for every
further limb of the root; over the range from degree~3 with small
coefficients to degree~14 with 400-bit ones it varies by a factor of ten. The
program estimates it from the degree, the largest coefficient and the height
bound before it looks at a prime --- the constants are
\texttt{RATPOINTS\_CHECK\_STEP}, \texttt{\_LIMB}, \texttt{\_CALL} and
\texttt{\_ROOT}, relative to \texttt{RATPOINTS\_CHECK\_REFERENCE}, the value
for a degree-6 curve with small coefficients --- and divides the two fractions
above by the ratio, so that a curve whose test is dearer gets more primes.
Only ratios enter, so a machine on which every test costs more needs no
change; \verb+-W+~$w$ overrides the estimate with $w$ rdtsc cycles
(Section~\ref{fd}), and \verb+-v+
prints what the estimate gives for the curve in hand.

\subsection{Correcting the choice while the program runs}\label{sec:adapt}

The rules above use the densities alone, and they miss one thing by
construction: if $a/b$ gives a rational point, then so does $ka/kb$ for
every~$k \ge 1$, and this pair therefore passes
every prime test there is. A~floor of survivors therefore outlives any amount
of sieving, and only the test for common factors removes it. The program can
see the floor: it counts the bit arrays that survive
the first stage, the survivors of the second, and the candidates that pass
the test for common factors, and fits
\[ S(n) = \text{floor} + \text{chance}\cdot R(n) \]
to the counts, with $R(n)$ the product of the $n$ smallest densities. It
revises its choice after the first million words of numerators and then
whenever twice as many have been swept. By default (\verb+-A 1+) it revises
only the third stage, where the number of survivors a denominator brings
becomes a count rather than an estimate; \verb+-A 2+ reconsiders
\texttt{sp2} as well, and \verb+-A 0+ switches the revision off.
Nothing that was fixed with \verb+-n+, \verb+-N+ or \verb+-P+ is revised.

\subsection{Retuning for your machine}\label{sec:retune}

The five constants --- the threshold, the offset, the run length at which a
second-stage prime pays for setting itself up, what building a table row
costs, and the third stage's cost per denominator --- were fitted on the
author's laptop by minimising the sum of the times of \texttt{make test1}
and \texttt{make test1many}, which between them cover the two regimes that
occur in practice, random curves and curves with many rational points. To
repeat that on your machine, run
\begin{verbatim}
> make tune
\end{verbatim}
It sweeps the five constants, in that order, over both test sets and writes
the values it finds to~\texttt{tuning.mk}, which the \texttt{Makefile}
includes, so that a subsequent \texttt{make all} builds with them. No source
file is touched, and deleting \texttt{tuning.mk} restores the values compiled
into~\texttt{ratpoints.h}. The file records the configuration it was
measured for and is ignored, with a warning, if you later build with another
register width or~\texttt{PRIME\_SIZE}.

\texttt{make tune} is deliberately not part of~\texttt{make all}: it takes several minutes
(about three for the 256-bit build on the author's machine, four for
128~bits, six for the plain 64-bit one and eight for the 512-bit code
emulated through narrower operations), wants an otherwise idle machine, and
must not be run under~\texttt{make -j}. Timing this reliably is the hard part.
The cost surface is flat --- anything within a factor of two of a good
threshold costs under $3\%$ --- whereas a laptop under sustained load slows by
a quarter as it warms up, so a straight sweep would measure the clock rather
than the settings. Each candidate is therefore timed back to back with the
current settings and only the ratio is kept, the median over several rounds
is used, and the current settings are themselves among the candidates, so
that the ratio they measure against themselves shows how much noise is left.
If that is more than $2\%$ the run is declared inconclusive, and nothing is
written unless the winner is better by more than $3\%$: a noisy run leaves a
good default alone. Two consequences are worth knowing in advance. A run that
reports that nothing beat the current settings has \emph{succeeded}; on a
build that is already tuned that is the right answer. And a configuration
whose tests take long enough for the machine to throttle within a single
round may report itself inconclusive however often it is run; raising
\texttt{ROUNDS}, or using a machine that does not throttle, is then the
remedy. The environment variables \texttt{ROUNDS} (how many times each
setting is measured, default~$3$), \texttt{WARMUP} (seconds of load before
anything is counted, $20$), \texttt{MARGIN} ($0.97$) and \texttt{NOISE}
($0.02$) adjust this, and \texttt{R\_VALUES}, \texttt{E\_VALUES},
\texttt{U\_VALUES}, \texttt{C\_VALUES} and \texttt{Q\_VALUES} give the
ladders of candidates for the five constants.

\begin{verbatim}
> make tunehigh
\end{verbatim}
measures on \texttt{make testhigh} and \texttt{make testhighmany} instead,
at a height bound of $2\cdot10^5$, where the two sieving stages are nine
tenths of the run; at the $16383$ of \texttt{test1} two fifths of the time on
the random curves goes into work other than those two stages, and the
constants have no effect on most of that.
Use it if the runs that matter to you are long ones. One timing there is two
minutes rather than three seconds, so it does not sweep the whole ladder: it
starts from the settings in force --- which is to say from what \texttt{make
tune} found, since both read and write the same \texttt{tuning.mk} --- and
asks only whether a factor of two either way in the threshold, in the run
length, in the table cost or in the third stage's cost per denominator, or
two more or fewer second-stage primes, is better (\texttt{R\_FACTORS},
\texttt{E\_DELTAS}, \texttt{U\_FACTORS}, \texttt{C\_FACTORS} and
\texttt{Q\_FACTORS} set that neighbourhood). Expect a couple of hours at the
default three rounds, and run it again if it moved a value. The intended
order is \texttt{make tune} first and \texttt{make tunehigh} afterwards;
only a pair of tunings at run lengths a hundred times apart can separate the
offset from the run length it is scaled by.

Should you prefer to do it by hand, the constants can be given for a single
run with \verb+-r+, \verb+-R+, \verb+-U+, \verb+-C+ and \verb+-Q+, which
\texttt{ratpoints} and \texttt{rptest} both accept, and hard-coded by
changing the definitions in \texttt{ratpoints.h} (followed by \texttt{make
distclean all test}). \texttt{rptest -z -T} prints the CPU time of a run and
nothing else, which is what \texttt{tune.sh} uses; \verb+-h+~$H$ changes its
height bound, \verb+-m+~$m$ repeats the computation $m$~times, and
\texttt{testdata.h} can be replaced by a collection of your own curves.
Programs using the library set the same constants through the fields of
\texttt{ratpoints\_args} (Section~\ref{sec:library}).


\section{Change log}\label{sec:changelog}

\textbf{Version 2.0} was released January 9, 2008.

\textbf{Version 2.0.1} was released July 7, 2008. It fixes a bug that prevented the
'\texttt{-1}' option to work properly.

\textbf{Version 2.1} was released March 9, 2009. It makes use of the SSE instructions,
so that the sieve can work on 128~bits in parallel (instead of on 64 or~32).
On my laptop (with an Intel Core2 processor), \texttt{make test} runs about
25\% faster than before.

\textbf{Version 2.1.1} was released April 14, 2009. It fixes a bug that in some
cases prevented an early abort when \texttt{*quit} was set by the callback function.
Thanks to \textbf{Robert Miller} for the bug report and the fix. In addition, it is now
checked that \verb+__WORDSIZE == 64+ before SSE instructions are used
(in \texttt{rp-private.h}), since the program then assumes that a bit array consists
of two \texttt{unsigned long}s. In this context, \verb+__SSE2__+ has been replaced
by a new macro \verb+USE_SSE+. A further fix eliminates unnecessary copying
of sieving information when SSE instructions are not used (introduced in version~2.1).
This was almost always harmless, but could have resulted in memory corruption
in the extreme case that all the information had to be computed for all the
primes.

\textbf{Version 2.1.2} was released May 27, 2009. It fixes some memory leaks.
Thanks to \textbf{Robert Miller}.

\textbf{Version 2.1.3} was released September 21, 2009. The library function
\verb+find_points+ should now work without any bound on the degree of the
polynomial. The degree bound for the \texttt{ratpoints} executable is now set
to~100 by default (specified by \verb+RATPOINTS_MAX_DEGREE+, used to be~10).

A bug in \texttt{rptest.c} (pointed out to me by \textbf{Giovanni Mascellani}
and \textbf{Randall Rathbun}) that could lead to a segmentation fault was fixed on
March 10, 2011.

\textbf{Version 2.2} was released January 9, 2022. The main change is that
\texttt{ratpoints} can now also work with 256-bit (and, in principle, 512-bit)
registers when the CPU has the relevant capabilities. Thanks to
\textbf{Bill Allombert} for doing a first implementation of this. In addition,
the code has been cleaned up to some extent (e.g., the variants according
to register sizes are now completely dealt with through macros defined
in \texttt{rp-private.h}), and some more comments have been added.
There are also some further optimizations; for example, the first sieving
stage now uses many of the SSE/AVX registers in parallel, and on some
architectures, the code uses a faster implementation for the test whether
a bit array is zero. On my current laptop, the new code (using 256-bit words)
is about twice as fast as the old one (2.1.3 using 128-bit words).

\textbf{Version 2.2.1} was released January 18, 2022. This fixes a small
bug that was introduced in the previous version, which could lead to
the program missing points at the upper end of the search interval
for the numerators. Thanks to \textbf{Bill Allombert}.

\textbf{Version 2.2.2} was released May 15, 2023. This fixes a small
bug that led to the polynomials being reversed incorrectly when the
degree is odd. Thanks to \textbf{Nicolas Mascot} and \textbf{Bill Allombert}.

Work on the code starting with version~2.2.3 was done with the help
of \textbf{Claude Code} (by Anthropic), using various models
(Opus~5, Fable~5, and Fable~5.1).

\textbf{Version 2.2.3} was released September 6, 2026. It fixes the code
path using 512-bit registers, which had never been exercised. Two problems
were found. The first one would have made
the program crash. The second problem was that the \texttt{TEST} macro, which is
used in the innermost loop of the second stage, still used the generic
fall-back version, leading to a degradation in performance;
it now uses \texttt{vptestmq}, in analogy with the 256-bit version.

Since I have no CPU with AVX512F available, this was tested by compiling
with \verb+-DUSE_AVX512+ but without \verb+-mavx512f+, so that
\texttt{gcc} expresses the 64-byte vectors through 256-bit operations.
This exercises the whole code path. What this cannot check is the
genuine 512-bit instructions. In the meantime, \textbf{Drew Sutherland} has
run the two variants against each other on an AMD Ryzen~9 9950X3D (Zen~5
architecture), for a curve with many rational points and a height bound
of~$10^6$; the 512-bit version was faster by 13 to 14\%.

\textbf{Version 2.2.4} was released September 19, 2026. It fixes bugs that a
review of the code and running the new test suite turned up:
\begin{itemize}
  \item A polynomial of degree~1 --- given as such,
        or reached from $y^2 = c_2x^2 + c_1x$ by the reversal, or from a zero
        leading coefficient --- made the Sturm sequence code read past the end of
        its arrays and crash; this affects library users as well.
  \item When the positivity region of~$f$ did not meet the search domain, the program
        returned before it had looked at the points at infinity, so that
        \texttt{ratpoints '-1000000000 0 0 1' 45} printed nothing where \verb+-s+
        printed $(1:0:0)$, and after a reversal the lost point is an affine point
        of the given curve.
  \item When no residue class of the denominator modulo~16
        admitted any numerator, the array of numerator patterns was read before it
        had been filled; by default the answer was still right, but with
        \verb+-j -F 0+ the program printed nondeterministic points outside the
        height bound after a run several thousand times longer than necessary.
  \item Several input fields of \texttt{ratpoints\_args} were used as scratch
        space --- the number of forbidden divisors found was written into
        \texttt{max\_forbidden}, the intersected search intervals into
        \texttt{num\_inter}, normalised values into \texttt{sp1}, \texttt{sp2}
        and others --- so that a program which fills the structure once and then
        calls \texttt{find\_points\_work} for curve after curve, as the section on
        the library suggests, had the forbidden-divisor test switched off after the
        first curve without forbidden divisors (a factor of two in running time)
        and the first curve's positivity region intersected into the search domains
        of the others, which loses points. The program no longer writes either,
        and the manual now says exactly which fields it does set;
        \texttt{make test1once} runs the test curves with the structure filled
        once, and \texttt{make test3} runs one or more invocations of
        \texttt{ratpoints} per bug against a reference.
  \item The library tested whether the leading coefficient is a square before
        it checked the coefficient pointer for \texttt{NULL}, so that a caller
        passing \texttt{NULL} crashed instead of getting the error code back.
  \item An empty coefficient string made the program report
        ``Bug no.~1 --- please report!'' instead of refusing the input.
  \item At the very top of the range of a \texttt{long}, three things went
        wrong: a search interval that reached the last few hundred numerators
        below~$2^{63}$ was dropped without a word, two computations on the
        position of its last bit having overflowed; the loop over every
        denominator never ended when the denominator bound was $2^{63}-1$,
        since the increment wraps around; and the loops over the square
        denominators formed $b^2$ and $d\,b^2$, which overflow within
        $2\cdot10^9$ of it.
  \item The test on the valuation of a denominator $d\,b^2$ at a prime
        dividing the leading coefficient applied \texttt{abs}, the
        \texttt{int} function, to a \texttt{long}, so that for denominators
        above~$2^{31}$ it judged the low 32~bits only and could exclude a
        denominator that carries a point: a point with denominator
        $3\cdot 59827^2$ was missed.
  \item Not a bug, but wasted time: the loops over the square denominators
        started at~1 whatever the lower bound \verb+-dl+ said; they now start
        at the first square in the range.
  \item Also fixed: four left shifts of negative values
        (undefined in C, harmless with gcc so far), two arrays that could be
        declared with length zero under \verb+-n 0 -N 0+, an overflow in the
        enumeration of the divisors of the leading coefficient for denominator
        bounds above about $4\cdot10^{15}$, and two constants in
        \texttt{ratpoints.h} that could not be overridden from the compiler
        command line. The library is binary compatible with 2.2.3.
\end{itemize}

\textbf{Version 3.0.0} was released September 21, 2026. Its theme is that
the program works out for itself what it used to be told, and the sieve has
been reworked from the set-up to the exact test. Measured against
version~2.2.4 on the same curves and the same laptop (cycles pinned to one
core, medians of three interleaved rounds; 2.2.4 at its own defaults, which
sieve with the primes below $128$ where this version uses those below $256$;
both find the same points on every suite):
\begin{center}
\begin{tabular}{lrrrr}
suite & 2.2.4, Gcycles & 3.0.0 & ratio & faster by \\
\hline
\texttt{make test1} (random curves, $16383$)      &   8.45 &   2.89 & 0.34 &  2.9 \\
\texttt{make test1many} (point-rich, $16383$)     &  11.27 &   4.80 & 0.43 &  2.3 \\
\texttt{make testhigh} (random, $2\cdot 10^5$)    & 412.7  & 190.8  & 0.46 &  2.2 \\
\texttt{make testhighmany} (point-rich, $2\cdot 10^5$) & 764.4 & 138.3 & 0.18 & 5.5 \\
random curves at height $4000$                  &   3.06 &   0.54 & 0.18 &  5.6 \\
random curves at height $1000$                  &   2.06 &   0.21 & 0.10 & 10.1 \\
point-rich curves at height $1000$              &   0.44 &   0.15 & 0.34 &  3.0 \\
random curves at height $200$                   &   0.96 &   0.12 & 0.13 &  8.0 \\
\end{tabular}
\end{center}
Instructions fell to between $0.16$ and $0.47$ of 2.2.4's and mispredicted
branches to between $0.03$ and $0.40$: at the small height bounds 2.2.4 was
bound by its branches, and cycles fall two to three times more than
instructions there. The changes, in outline; the sections named have the
details, and Appendix~\ref{app:notes} the measurements behind them.
\begin{itemize}
  \item \emph{The sieving parameters are chosen from the curve.} How many
        moduli each of the first two stages uses, and which ones, follows
        from the densities of admissible residues, from the length of the
        run and from what each modulus costs, with the moduli ranked by
        cost per information rather than by information alone
        (Section~\ref{sec:spchoice}); \verb+-n+ and \verb+-N+ still fix the
        numbers. When the small primes say nothing, as on a curve with very
        many rational points, the program looks at further primes of its
        own accord, and the default table holds the primes below~$256$
        (\texttt{PRIME\_SIZE}~8, formerly~7) to give it something to reach
        for. The rules rest on five machine constants, set with
        \verb+-r+, \verb+-R+, \verb+-U+, \verb+-C+ and \verb+-Q+ and
        measured by \texttt{make tune} and \texttt{make tunehigh}
        (Section~\ref{sec:retune}), which write them to \texttt{tuning.mk}.
  \item \emph{A third sieving stage} tests the survivors of the two
        bit-array stages one at a time against further primes, without
        sieving tables, so that it can use primes far up the list; how many
        is decided per curve (Section~\ref{sec:sp3choice}, option~\verb+-P+)
        from an estimate of what one exact test costs on that curve
        (\verb+-W+). While it runs, the program counts what the sieve is
        finding and revises the choice (Section~\ref{sec:adapt},
        option~\verb+-A+).
  \item \emph{Composite moduli.} Odd composite numbers below~$64$ sieve alongside
        the primes (the bound is given by~\texttt{RATPOINTS\_COMPOSITE\_MAX}
        in~\texttt{ratpoints.h}). A prime at which $f$ is a square at every
        residue is no longer discarded when denominators divisible by it
        occur: it still says that numerator and denominator are not both
        divisible by it, and it says so at no cost as a factor of a
        composite modulus (Section~\ref{sec:spchoice}).
  \item \emph{2-adic information modulo~$64$} instead of~$16$, the even
        denominators decided exactly through the reversed polynomial, and
        the bit arrays of each denominator class packed with the largest
        stride $2^k$ its numerators allow rather than by parity alone
        (Section~\ref{impl}).
  \item \emph{Denominators.} The test at a prime dividing the leading
        coefficient, derived from the Newton polygon of~$F$, now runs;
        the forbidden-divisor bit
        arrays cover every excluded prime up to the square root of the
        height bound, and the default of \verb+-F+ is~$64$; the Jacobi
        symbol is a product of table-driven Legendre symbols; the tests on
        $b \bmod 64$ are applied to 64~denominators at a time; the residues
        of~$b$ come from precomputed reciprocals, and the gcd that tests a
        candidate for lowest terms is branchless.
  \item \emph{The sieving loops.} The first stage \texttt{AND}s the 2-adic
        pattern in with its first modulus instead of pre-filling the bit
        arrays, and sieves the bit arrays left over after the last whole
        chunk of sixteen in legs of 8, 4, 2 and~1 registers instead of
        padding the range; the second stage steps over the empty bit arrays
        into a sentinel, reaches its table rows through a precomputed
        reciprocal and extracts the survivors by counting trailing zeros.
        The sieving tables are built four times faster, and the memory
        reserved for them follows the number of primes in use.
  \item \emph{Build and portability.} A stamp file makes a change of
        \texttt{CCFLAGS1} or \texttt{PRIME\_SIZE} rebuild what depends on it;
        \texttt{CCFLAGS128} selects \verb+USE_AVX128+, which needs only SSE2
        and is faster than the \verb+USE_SSE+ variant; the code assumes a
        64-bit \texttt{long} and refuses to compile without one, so
        version~2.2.4 is the last one for 32-bit machines.
  \item \emph{Library.} \texttt{ratpoints\_args} has gained fields, so
        programs using the library must be recompiled; the input fields
        come back as they went in, except \texttt{cof} and \texttt{degree},
        which describe the polynomial the program worked with, and
        \texttt{num\_inter} and \texttt{domain}, which describe the region
        searched; and
        \texttt{sp1\_used}, \texttt{sp2\_used} and \texttt{sp3\_used} report
        what was used (Section~\ref{sec:library}).
  \item \emph{Tests.} \texttt{make test} runs, besides the random curves,
        the point-rich ones (\texttt{test1many}), a hundred curves of
        degrees 3, 4, 7 and~8 (\texttt{testdegrees}), the argument structure
        filled once for all curves (\texttt{test1once}), one invocation per
        bug fixed (\texttt{test3}), a suite that exercises nearly every line of
        the code against a reference checked by brute force (\texttt{test4})
        and the library interface (\texttt{testapi}), and it fails when a
        test fails; \texttt{testhigh}, \texttt{testhighmany},
        \texttt{test4configs} and \texttt{coverage} are run separately
        (Section~\ref{sec:tests}).
\end{itemize}

\textbf{Version 3.1.0} was released September 24, 2026. The sieve can run on
several threads.
\begin{itemize}
  \item \emph{Threads.} The option \verb+-t+~$n$ and the field
        \texttt{num\_threads} of \texttt{ratpoints\_args} run the sieve
        and the exact tests on $n$ threads besides the calling one, which
        hands out the work and delivers the points
        (Sections~\ref{sec:threads} and~\ref{sec:library}). The loop over
        the denominators is cut into blocks whose length depends on the
        search alone, the threads take the blocks in order, and the points
        of a block are delivered after those of every earlier block: what
        the program prints, what the callback receives and in which order,
        and where \texttt{*quit} takes effect do not depend on the number
        of threads, and neither do the sieving parameters chosen and their
        revision during the run, which follow the same schedule of blocks.
        The sieving tables are built once and shared. POSIX threads;
        \texttt{make THREADS=0} builds without them, and at most $33$
        threads sieve unless the library is built for more
        (Section~\ref{sec:build}). A search that cannot
        keep the points of a block fails with \texttt{RATPOINTS\_NO\_MEMORY}
        rather than deliver fewer points than it found. On a 16-core
        machine a search at height bound $4\cdot10^5$ runs $9.4$ times
        faster on 16 threads than on one and the long test suites $5.7$
        to $7.9$ times faster. Programs using the library must
        be recompiled, since \texttt{ratpoints\_args} has gained fields.
  \item \emph{The 512-bit build} has been run on AVX-512 hardware, where the
        test suite passes; see Section~\ref{sec:build} for when it makes sense
        to use it.
  \item \emph{Tests.} \texttt{make testthreads} runs the suites with several
        threads and with the blocks of the denominator loop fixed at other
        lengths, and \texttt{make tsan} runs them under ThreadSanitizer
        (Section~\ref{sec:tests}).
\end{itemize}


\appendix

\section{Notes on the optimizations}\label{app:notes}

This appendix was written by Claude Code, the AI assistant (by Anthropic)
with whose help this version was developed (Section~\ref{sec:changelog}); it
consolidates the notes kept during that work. It keeps what a reader would
need in order to change the program without repeating the experiments: how
the measurements were made and what they can and cannot resolve, where the
time goes, why curves with many rational points behave so differently from
random ones, the cost model with the origin of each constant, what each
change was worth, what was tried and did not pay, and how much of all this
depends on the machine. Every figure was taken on one laptop --- an Intel
i7-1355U, the code pinned to one performance core, the 256-bit build,
gcc~14 --- and is quoted as a ratio of cycle counts, the program with the
change in question over the program without it;
Section~\ref{app:machine} says what would come out differently elsewhere.

\subsection{How the measurements were made}\label{app:measure}

The effects in question are mostly between one and ten per cent, and the
machine is noisier than that. Under sustained load it slows by about a
quarter as it warms up, over tens of seconds, so wall-clock times taken
minutes apart mean nothing; core cycles counted by \texttt{perf} on one
pinned core do better, but the drift affects them too. Every comparison was
therefore made in pairs: the candidate and the baseline run back to back,
interleaved, for several rounds; only the ratio of each round is kept, and
the median of the per-round ratios is quoted, not the ratio of the medians
(the two differ, because the drift is not symmetric between the two
positions). The baseline is run against itself as one of the candidates,
and its ratio measures the noise that is left; a run in which that exceeds
two per cent is not believed. \texttt{tune.sh} does all of this and refuses
to write anything a noisy run suggests.

Two floors remain below which whole-program timings cannot decide anything.
Where the compiler places a hot loop is worth up to ten per cent on this
machine: the same source rebuilt with \texttt{-falign-loops=32} and~\texttt{64}
gave ratios of $0.98$ and $1.00$ where the default placement gave $1.09$.
Whole-program comparisons were therefore made at three placements of the
code, and a change under about two per cent at a height bound of
$2\cdot10^5$ --- where half the run is the innermost loop of the first stage,
which is exactly what those flags move --- cannot be settled that way at
all. The other method, putting both variants in one binary behind a switch
read once at start-up, removes the placement effect entirely and is the
better evidence; but it is blind to a data-layout cost. A change that grows
a structure read in the innermost loop by eight bytes shows nothing within
one binary and gives back half of a $2.8\%$ gain across builds. Both methods
were used when a change touched such a structure.

Instruction counts are exact, repeatable to $10^{-8}$, and the right tool for
attributing a set-up cost to a function; as a proxy for time they misled
three times. \texttt{-funswitch-loops} removed one to three per cent of the
instructions and gained nothing in cycles at any placement; hoisting the
per-curve constants of the sieve entries out of the per-denominator loop
removed $2.4\%$ of the instructions and $0.1\%$ of the cycles, the loop being
bound by its chain of dependent loads; and composite moduli above~$64$
removed $22\%$ of the instructions of \texttt{make testhigh} while the cycles
fell less than with the bound in place and \texttt{make testhighmany} lost
$7\%$, because the rows of a large modulus do not stay in the first-level
cache. Cycles decide; instructions explain.

The test suites see different things. At the height bound $16383$ of
\texttt{make test1} and \texttt{test1many} two fifths of the time on the
random curves is not sieving at all but choosing the primes, building the
tables and setting up each denominator; at the $2\cdot10^5$ of
\texttt{testhigh} and \texttt{testhighmany} the two sieving stages are nine
tenths of the run. A change to the sieving loops is therefore judged on the
latter, a change to the set-up on the former or at height bounds of $200$
to~$4000$, and the two regimes of Section~\ref{app:regimes} are always
reported apart, since they often disagree and the disagreement is usually
the finding. One more thing a correctness suite cannot see: a change that
removes an initialisation reproduces every reference output, because the
values it stops writing were about to be overwritten --- and a debug or
instrumented build may still read them. The \texttt{-DDEBUG} and the
\texttt{-DRP\_PHASE\_TIMING} builds are run under \texttt{valgrind} after
such a change.

\subsection{Where the time goes}\label{app:profile}

The sieving tables are built once for each pair of a modulus and a class
of denominators modulo it and kept, so their cost is bounded by the moduli
and does not grow with the height bound; everything else grows at least
linearly with it. The shares therefore depend on the height bound above
all. On the random test curves:
\begin{center}\small
\begin{tabular}{lrr}
 & height $16383$ & height $2\cdot10^5$ \\
\hline
the two sieving stages, with the exact checks & $57\%$ & $90\%$ \\
the sieving tables & $4\%$ & $0.3\%$ \\
the set-up per denominator & $8.5\%$ & $1.3\%$ \\
outside the sieve altogether: examining and choosing the moduli & $16\%$ & $3\%$ \\
\end{tabular}
\end{center}
(the rest is the work of \texttt{sift()} around the stages). Within the
sieve the first stage is the larger part, of the order of three fifths of a
long run on random curves and three quarters on point-rich ones; the second
stage with the exact checks is most of the rest. The exact check itself is
under one per cent of every suite, and the integer divider is busy for less
than $0.6\%$ of the cycles: neither is where the time goes
(Section~\ref{sec:remarks}). At a height bound of~$200$ the shape is
different again: examining the thirty primes is a quarter of the
instructions, and the sieving tables are a twentieth of the run --- they
would be a fifth if the first-stage rule did not price them
(Section~\ref{sec:spchoice}).

The first stage's inner loop is close to what the machine can do: sixteen
fused load-and-\texttt{AND} instructions per modulus per chunk and about
eight instructions of overhead, at $1.8$ instructions per cycle, bound by
the load throughput and the capacity of the first-level cache rather than
by arithmetic. Its working set at the parameters a random curve uses is
some $50$~KB at $256$~bits, which is about the size of that cache on the
development machine; the marginal cost of a first-stage prime is flat while
the set fits and rises past it. An \texttt{AND} costs one cycle while the
rows of the first-stage moduli together occupy under some $25$~KB, and about
$0.8\%$ more for every further KB. The mispredicted branches, of which
version 2.2.4 had many, scale with the survivors and not with the
\texttt{AND}s.

\subsection{The two regimes}\label{app:regimes}

Almost every tuning question in the program turns on the proportion of
numerators that survive the first sieving stage, and two populations of
curves differ very much in what it takes to bring that proportion down. The
rule of Section~\ref{sec:spchoice} lets about one numerator in three
thousand through at a height bound of~$16383$, and one in five or six
thousand at~$2\cdot10^5$, on either kind of curve; but on a random curve it
gets there with about ten moduli, while on a curve with very many rational
points $f$ is a square modulo every residue for the smallest primes, so
those primes carry next to no information --- only that numerator and
denominator are not both divisible by them --- and the rule needs
fifteen to twenty moduli. Of the $98$ point-rich test curves, $70$ make the
program look beyond the thirty primes it examines first. That is
why the program looks at further primes when it needs to, and why the table
holds the primes below~$256$: it is worth a factor of $1.3$ to~$2$ on those
curves at a large height bound and costs nothing on any other. Larger
primes are what these curves need: the proportion of admissible residues
tends to one half as the prime grows however many points the curve has,
while for the small primes it is close to one.

The second thing that separates the regimes is a floor of survivors that
no sieving removes: if $a/b$ is a rational point then so is $ka/kb$ for
every~$k$, with the same value of~$f$, and it passes every prime test. There
are about $H/b$ of them for each point of denominator~$b$, so they grow
with the height bound while the survivors that die by chance grow with its
square; only the test for common factors removes them, which is why that
test runs before the third stage. The floor is what the densities alone
cannot predict and what the correction of Section~\ref{sec:adapt} measures,
and it is why the number of candidates per denominator that survive the
second stage differs by a factor of ten between the suites --- $0.5$ and
$0.33$ on the random curves at $16383$ and at $2\cdot10^5$, $2.7$ and $4.1$
on the point-rich ones --- so that no fixed number of third-stage primes can
serve and the rule there is a stopping rule.

The two regimes want the same constants, which was not obvious and was
measured twice: \texttt{make tunehigh}, which tunes at $2\cdot10^5$, keeps
every value that \texttt{make tune} finds at $16383$, all within $0.7\%$
except the threshold, which is sharper --- doubling it costs $5.6\%$ at
$2\cdot10^5$, and two and a half times its value costs $8.6\%$ at $16383$.
One threshold also serves every height bound from $200$ to $2\cdot10^5$.
That depends on the estimates behind the rules --- the words the run will
sweep, the denominators it will visit, the bits set per word --- being right
together (Section~\ref{app:model}): with one of them off, every height
bound wants a threshold of its own.

\subsection{The cost model and where the constants came from}\label{app:model}

The rules of Section~\ref{sec:tuning} compare costs, and every cost in
\texttt{ratpoints.h} is in one unit: what one first-stage \texttt{AND} costs
per 64-bit word of numerators, a quarter of one \texttt{AND} on a 256-bit
array, about $0.26$ core cycles on the development machine. The constants
that are not fitted were measured with the instrumented build
(\texttt{-DRP\_PHASE\_TIMING -DRP\_PHASE\_COUNTS}), dividing each part's
cycles by the number of times it ran. The instrumentation is not free ---
the counters alone can be nearly a fifth of an instrumented run --- so shares are
taken from a timing-only build. Table~\ref{tab:constants} lists the
constants: what each one is, how its value was found, and how much a
different value costs.

\begin{table}[tbp]\footnotesize\setlength{\tabcolsep}{3pt}%
\renewcommand{\arraystretch}{1.1}
\begin{tabular}{@{}P{4.0cm}P{1.7cm}P{5.6cm}P{5.6cm}@{}}
constant (\texttt{RATPOINTS\_\dots}) & value & what it is; how the value was found & how flat \\
\hline
\texttt{SURVIVORS\_PER\_WORD} & $0.003$ &
  the threshold at which the first stage stops taking moduli; fitted by
  \texttt{make tune} at $16383$ &
  $0.0045$ costs $2.7\%$ and $0.0075$ costs $8.6\%$; $0.0015$ to $0.002$
  cost nothing; the same at $2\cdot10^5$ \\[3pt]
\texttt{SP2\_EXTRA}, \texttt{SP2\_U0} & $11$, $1.6\cdot10^6$ &
  the second-stage offset of an arbitrarily long run, and the run length at
  which a prime's set-up costs as much as its sieving; fitted to four
  measured optima --- 3, 5, 9 and 12 second-stage primes at $U =
  5.9\cdot10^5$, $5.1\cdot10^6$, $8.2\cdot10^7$ and $7.1\cdot10^8$ words ---
  which the pair reproduces as 3, 8, 11, 11 &
  a flat offset (\verb+-U 0+) costs about $5\%$ at one height bound or the
  other; a flat~$5$ costs $12\%$ of \texttt{testhighmany}; only
  \texttt{tune} and \texttt{tunehigh} together separate the two \\[3pt]
\texttt{COST\_TABLE} & $38$ &
  one row of a sieving table; measured at 18 to 30 &
  flat from 20 to 70: $0.9515$ at 38 against $1$ at 0, $0.9619$ at 140 \\[3pt]
\texttt{COST\_BP}, \texttt{COST\_SETUP} & $8$, $30$ &
  a denominator's entry in the list of reduced denominators, and the fill
  of a sieve entry per denominator; measured &
  8 against 24 within $0.3\%$ \\[3pt]
\texttt{COST\_CALL}, \texttt{COST\_LINE} & $50$, $4.5$ &
  a call of the sieve per modulus (12 to 13 cycles) and the fetch of a table
  row per denominator (1.2 cycles per cache line); measured &
  not swept singly; fitted as one group with the threshold \\[3pt]
\texttt{COST\_PHASE2}, \texttt{COST\_SURVIVOR} & $110$, $1400$ &
  one second-stage \texttt{AND}, and one survivor of the second stage; the
  first measured at 155 (random) and 61 (point-rich) at $2\cdot10^5$ &
  60, 100 and 220 within about $2\%$ of 110 \\[3pt]
\texttt{SP3\_PER\_SURVIVOR}, \texttt{SP3\_PER\_DENOM} & $0.055$, $0.013$ &
  the third stage's costs per survivor tested and per denominator, as
  fractions of one exact check; measured: 250 to 280 cycles for a check
  beyond the first of a denominator, some 190 once per denominator, 3.5 per
  step carrying $b$ past a prime &
  the per-denominator value from $0.003$ to $0.05$ within $1.2\%$; a fixed
  \verb+-P+ never beats the rule by more than $0.3\%$ \\[3pt]
\texttt{CHECK\_STEP}, \texttt{\_LIMB}, \texttt{\_CALL}, \texttt{\_ROOT} &
  $29$, $8$, $94$, $170$ &
  the exact check, in rdtsc cycles: a Horner step, a further limb, the
  call, a further limb of the root; from \texttt{make bench\_check} (in the
  sieve the check costs $44 + 1.44$ times the microbenchmark) &
  only ratios between curves enter; the estimate agrees with the sieve to
  $6\%$ over degrees 3 to 12 \\[3pt]
\texttt{CHECK\_REFERENCE} & $306$ &
  the same for a degree-6 curve with small coefficients, which is what
  \verb+-W 306+ assumes for every curve &
  scaling every curve's value alike from 306 to 1000 moves no suite by
  $0.5\%$ \\
\end{tabular}
\caption{The constants of the cost model: their values, where they came
from, and how much a different value costs.}\label{tab:constants}
\end{table}

Two properties of the model are worth more than any single value. First,
the rules sit at their own optimum, so a constant off by a factor of two
costs under three per cent, and a better estimate of one moves the boundary
by about one prime, which is by definition worth about what it costs; that
is why the per-curve estimate of the exact check's cost changes no running
time by half a per cent and is there so that the constants mean what they
say. Second, the estimates the rules feed on --- the words the run will
sweep, the denominators it will visit, the bits set per word, the density
of a modulus --- and the fitted constants form one system: the constants are
fitted with the estimates in place and absorb whatever error the estimates
have, so a more accurate estimate on its own makes the program slower, and
has to be judged after a refit of the constants (\texttt{make tune}), not
before. As they stand, the predicted words per denominator are within
$1.5\%$ of the count at every height bound from $200$ to $2\cdot10^5$, and
the automatic choice of the first-stage moduli beats or ties the best fixed
\verb+-n+ in every one of eighteen combinations of suite, height bound and
modulus cap. A finer point of the same kind is the density of a prime: with
constants that had been fitted to an approximate density (the class of
denominators divisible by~$p$ counted as admitting every numerator), the
exact one made the point-rich curves $5\%$ slower at $2\cdot10^5$, and it
became the better choice only once the constants were refitted with it. A
fitted model absorbs its approximations, and ``more exact'' can be worse
until the fit is repeated.

\subsection{What each change was worth}\label{app:worth}

Table~\ref{tab:worth} gives, for each change, the cycles with it over the
cycles without it on the four suites --- random and point-rich curves at
the height bound $16383$ (\texttt{make test1} and \texttt{test1many}) and at
$2\cdot10^5$ (\texttt{testhigh} and \texttt{testhighmany}) --- and, where it
matters, at other height bounds or on the suite of curves of degrees 3, 4,
7 and~8 (\texttt{make testdegrees}). Each was measured as described above,
against the program as it stood when the change was made, so the rows do
not multiply out to the totals of Section~\ref{sec:changelog}. A dash is a
change within the noise; an empty cell was not measured.

\begin{table}[tbp]\footnotesize\setlength{\tabcolsep}{3pt}%
\renewcommand{\arraystretch}{1.1}
\begin{tabular}{@{}P{7.3cm}rrrrP{3.8cm}@{}}
 & \multicolumn{2}{c}{height $16383$} & \multicolumn{2}{c}{height $2\cdot10^5$} & \\
change & random & rich & random & rich & elsewhere \\
\hline
\texttt{sp1}, \texttt{sp2} chosen from the survival rate; the tight scan
  over the empty arrays & $0.86$ & $0.64$ & & & \\
further primes for starved curves; \texttt{PRIME\_SIZE}~8 & -- & -- & -- & --
  & $0.50$ to $0.75$ on the curves it applies to \\
survivors located by counting trailing zeros & $0.99$ & $0.95$ & $0.98$ & $0.98$ & \\
the table set-up: the residue test's branch becomes a shift, four rows at a
  time & $0.87$ & $0.98$ & -- & -- & the tables $3.7\times$ faster \\
the third sieving stage & -- & $0.98$ & $0.98$ & $0.93$ & \\
reductions by a multiply-high instead of a subtraction chain & $0.97$ & $0.97$ & $0.98$ & $0.98$ & \\
the run length known in advance: cost-aware ranking, the offset scaled, the
  survivors counted, the check's cost estimated & $0.94$ & -- & -- & $0.88$
  & $0.86$ on the degree suite \\
the first stage writes the arrays: no fill pass & $0.97$ & -- & $0.94$ & $0.97$
  & $0.95$ on the degree suite at $2\cdot10^5$ \\
the valuation test at primes dividing the leading coefficient & $0.92$ & $0.98$
  & $0.91$ & $0.98$ & $0.95$ on the degree suite \\
the scan into a sentinel; the second stage's row by the reciprocal & $0.97$
  & $0.97$ & $0.93$--$0.96$ & $0.93$--$0.96$ & \\
the denominator side: Legendre-symbol tables, 64 denominators a word,
  reciprocals, a branchless gcd, forbidden divisors up to $\sqrt H$, the
  third stage set up on demand & $0.76$ & $0.96$ & $0.94$ & $0.98$ & \\
no padding: the tail in legs of 8, 4, 2, 1 words & $0.90$ & $0.94$ & $0.98$
  & $0.98$--$1.00$ & $0.81$ at 1000 and 4000, $0.92$ at 200 \\
the 2-adic information modulo 64 & $0.97$ & $0.98$ & $0.97$ & $0.99$
  & $0.97$--$0.98$ at 200 to 4000 \\
the third stage looks further; even denominators at half width & -- & $0.99$
  & -- & -- & $4.6\times$ fewer exact checks on the point-rich curves at
  $2\cdot10^5$ \\
each denominator class packed with its 2-adic stride & $0.91$ & $0.99$
  & $0.875$ & $0.98$ & $0.97$ at 4000, $1.035$ at 1000, $1.05$ at 200 \\
composite moduli up to 64 & $0.92$ & -- & $0.89$ & $0.96$ & nothing below a
  few hundred \\
the first-stage rule prices the tables and the downstream cost & -- & -- & --
  & -- & $0.76$ at 200, $0.86$ at 1000, $0.99$ at 4000 \\
the estimates of the run (words, denominators, bits per word, densities)
  made accurate, the constants refitted & $0.995$ & $0.985$ & $0.993$
  & $0.981$ & $1.00$ to $1.03$ at 200 to 4000 \\
primes at which $f$ is a square at every residue kept for what they say
  about common factors, mostly as factors of composite moduli & -- & $0.985$
  & -- & -- & $0.987$ on the point-rich curves at 1000 \\
\end{tabular}
\caption{What each change was worth: cycles with it over cycles without
it, on the random and the point-rich curves at the height bounds $16383$
and $2\cdot10^5$; a dash is a change within the noise.}\label{tab:worth}
\end{table}

The degree suite gained most from the changes to the model because it is
the one set none of the constants had been fitted on: rules carry over
where fitted constants do not. The largest single gain on the random curves
at $16383$ came from the denominator side, where the Jacobi symbol computed
by a binary algorithm for every denominator had been $12\%$ of the run and
a third of its mispredicted branches; the largest at $2\cdot10^5$ from the
stride, which sweeps $15\%$ fewer bit arrays.

\subsection{What was tried and did not pay}\label{app:declined}

\begin{itemize}
  \item \emph{The second stage on narrower registers.} The guess was that a
        wider array holds more bits and so takes more primes to clear. The
        counters say the number of second-stage \texttt{AND} steps is the same
        to $1.5\%$ at 64, 128, 256 and 512 bits, and a hybrid that scans at
        the register width and sieves the surviving words one at a time
        (which is what \texttt{USE\_LONG\_IN\_PHASE\_2} selects) is a wash to
        $4\%$ slower: a narrow and a wide read of the same table come from
        one cache line, so there was never any traffic to save.
  \item \emph{Byte rotations in the table set-up.} Holding the $p$-bit pattern
        in eight shifted copies makes every row an unaligned load; but the
        rotation was never the expensive part --- three quarters of the set-up
        was the loop that builds the pattern, whose branch on a
        quadratic-residue test mispredicted half the time --- and the byte
        offset still steps through the same three-cycle chain. Worth $2.5\%$
        of a set-up that is a tenth of \texttt{make test1}, against eight
        scratch buffers and a byte-order assumption. The lesson that did pay:
        a loop carrying $x \mathrel{+}= d$; if $x \ge p$ then $x \mathrel{-}= p$
        costs three cycles a row whatever else it does, and the remedy is to
        run several of them at once.
  \item \emph{Sieving tables without the replication.} The tables hold each
        $p$-bit row \texttt{RBA\_PACK} times over so that every bit array
        read is aligned; dropping that cuts their footprint by nearly that
        factor and is legal, since vector loads need no alignment. Priced in
        isolation by shifting every table one word, the split loads cost
        $28\%$ of the first stage at 256 bits, against a cache cliff worth at
        most $17\%$ and only for curves that reach eighteen or more
        first-stage primes; worse with 512-bit registers, where seven offsets
        in eight split a line. Ordering the tables by prime, proposed for the
        same reason, is the case in effect already, since they are built
        lazily in the order of use; and the TLB is not a constraint.
  \item \emph{Testing several bit arrays at once in the scan} pays only when
        the proportion of surviving arrays varies widely, as it does when the
        number of first-stage primes is fixed; with the rate pinned at a few
        per cent by the rule, a group containing a survivor wastes its whole
        \texttt{OR}, and a group of one wins. Per-chunk survivor flags
        computed in the first stage were prototyped three ways: the branchy
        form mispredicted, the branchless forms spilled, since all sixteen
        vector registers hold live arrays. The sentinel was the answer.
  \item \emph{Compiler switches.} \texttt{-funswitch-loops} removes one to
        three per cent of the instructions and none of the cycles;
        profile-guided optimisation removed $6.5\%$ of the instructions and
        added two to eight per cent of cycles. Both declined.
  \item \emph{Composite moduli above 64}, and a footprint term in the cost
        model to price them: the products of two mid-sized primes save up to
        $45\%$ of the first stage's \texttt{AND}s on some random curves at
        $2\cdot10^5$ and $3.7\%$ of \texttt{make testhigh}, but a row of $220$
        bit arrays leaves the first-level cache and the point-rich suite
        loses $4.5$ to $7\%$. Telling the two apart needs a model of the row
        set with three or four cache constants; the bound stays. The 2-adic
        pattern modulo $256$ would remove one admissible class in twelve more
        for sixteen times the set-up.
  \item \emph{Reconsidering \texttt{sp2} from the counted survivors}
        (\verb+-A 2+) is a second answer to a question the offset scaled by
        the run length already answers, and whichever is applied later wins:
        without the offset it costs $5\%$ of \texttt{make testhighmany}, with
        it it adds nothing. Two answers to one question are an override, not
        a composition.
  \item \emph{The tie at $v_p(b) = 1$.} When $p$ divides the leading
        coefficient exactly once and not the next coefficient, a denominator
        with $v_p(b) = 1$ admits a single class of numerators modulo~$p$,
        which the sieve could use with $p - 1$ further tables per prime. The
        estimate assumed the prime is in the first stage; it usually is not,
        because the class it cannot sieve makes it look weak to the ranking.
        Refined to $0.4$ to $1\%$ of a run and not built.
  \item \emph{Blocking the prime loop for the cache} and \emph{a two-copy
        table layout} were proposed for cores with a 32~KB first-level cache
        or with 256-bit loads executed as two halves; on this machine's
        performance cores they lose, and for its efficiency cores an
        analysis of the loop predicted gains of $9\%$ and of up to $46\%$ of
        the first stage. Their sign depends on the cache size, so they belong
        behind a run-time gate and a measurement on such a machine, and were
        not built.
\end{itemize}

\subsection{How much of it is this machine}\label{app:machine}

Every gain this version added is general: it comes from the mathematics
(fewer numerators or denominators to look at, fewer rows to build) or from
the instruction stream (fewer instructions, fewer mispredicted branches, a
shorter dependency chain), and every out-of-order core rewards these in the
same direction. An algorithmic gain has the same size everywhere; an
instruction-stream gain scales with the core's throughput and its
misprediction penalty, which are within a factor of two of each other
across current cores. The machine lives in two places. The first is what
this version inherited: the register width, the sixteen registers of a
chunk (AVX-512 and ARM64 have thirty-two) and the replication of the tables,
the last two decided on a core with a 48~KB first-level cache where the
mainstream x86 cores of the last decade have 32. The second is the
constants: the \texttt{COST} values are ratios between different kinds of
work --- scalar integer code for a table row, a vector load for an
\texttt{AND}, branchy code for a survivor --- and their unit depends on the
register width, so that a 64-bit build has them $1.6$ times too small and a
128-bit build $1.2$ times; that is inside the flat basin, and \texttt{make
tune} covers the five that matter. Only the \texttt{TEST} macros are
intrinsics; the arithmetic on bit arrays is the compiler's vector
extension, so a port to another vector instruction set is one more branch
of some forty lines in \texttt{rp-private.h}.

\subsection{Traps for whoever changes the program}\label{app:traps}

\begin{itemize}
  \item Bit $j$ of \texttt{den\_bits} is the class $b \equiv j \bmod 64$, not
        $b - 1$, although the set-up shifts by $(b_{\text{low}} - 1) \bmod 64$:
        the loop shifts before it tests. Read the other way, it inverts the
        parity for a third of the curves and makes the run-length estimate
        $2.4$ times too small, which looks like a modelling error.
  \item Changing the number of moduli in the middle of a run is safe for
        correctness --- sieving only removes numerators that cannot be points
        --- but not for the bookkeeping: the list of moduli must hold every
        informative one, the list of reduced denominators must be sized for
        all of them, and the buffer of sieving tables must cover every prime
        looked at, since a prime the third stage held may be promoted and
        suddenly build a table. Missing the last is a segmentation fault.
  \item Anything added to the sieve entry structure is read in the innermost
        loop of the first stage; eight more bytes cost about $1.5\%$ of a
        long run.
  \item The sixteen-register ladder of the first stage stays hand-written:
        as an array of sixteen registers gcc spills it, $15\%$ slower. The
        tail legs need \texttt{\#pragma GCC unroll} to keep their smaller
        arrays in registers, and a variable that is in fact a constant must be
        written as one --- an offset written as \texttt{t \& \~{}15} instead
        of~\texttt{0} keeps eight index registers live and spills five per
        prime.
  \item \texttt{prime[]} holds the odd primes only, so \texttt{prime[0]} is~3.
  \item With \texttt{PRIME\_SIZE}~9 the mask that records which primes a
        modulus involves runs out of bits at the 64th prime, so beyond it the
        exclusion of moduli sharing a prime silently stops; the points are
        unaffected, the sieve merely does redundant work.
  \item The output of \verb+-x+, the survivors without the exact check,
        depends on which moduli were chosen and so on the tuning; only the
        points are invariant, and a test of the sieve must pin the primes.
\end{itemize}


\end{document}